\documentclass[onefignum,onetabnum]{siamart171218}

\usepackage{amsfonts}
\usepackage{amsmath}
\usepackage{amssymb}
\usepackage{mathtools}
\usepackage{graphicx}
\usepackage{epstopdf}
\usepackage[ruled,vlined,algo2e]{algorithm2e}
\usepackage{bbm}
\usepackage{lscape}
\usepackage{yfonts}
\usepackage{verbatim}
\usepackage{xcolor}
\usepackage{hyperref}
\usepackage{subcaption}

\ifpdf
  \DeclareGraphicsExtensions{.eps,.pdf,.png,.jpg}
\else
  \DeclareGraphicsExtensions{.eps}
\fi

\DeclareSymbolFont{lettersA}{U}{txmia}{m}{it}
\DeclareMathSymbol{\R}{\mathord}{lettersA}{"92}
\DeclareMathSymbol{\C}{\mathord}{lettersA}{"83}
\newcommand{\mfs}{\mathfrak{s}}
\newcommand{\mff}{\mathfrak{f}}

\newcommand{\bigO}[1]{\ensuremath{\mathop{}\mathopen{}\mathcal{O}\mathopen{}\left(#1\right)}}

\renewcommand{\sf}{\ensuremath{\{\textswab{s:},\textfrak{f}\}}}

\newcommand{\pl}{\ensuremath{\left(}}
\newcommand{\pr}{\ensuremath{\right)}}
\newcommand{\eff}{\ensuremath{\textnormal{eff}}}
\newcommand{\tol}{\ensuremath{\textsc{tol}}}

\definecolor{ForestGreen}{rgb}{0.18,0.70,0.18}

\newcommand{\editmade}[1]{#1}
\newcommand{\editblock}[1]{#1}

\allowdisplaybreaks

\newsiamremark{remark}{Remark}
\newsiamremark{hypothesis}{Hypothesis}
\crefname{hypothesis}{Hypothesis}{Hypotheses}
\newsiamthm{claim}{Claim}

\headers{Adaptive time step control for multirate infinitesimal methods}{A.~C.~Fish and D.~R.~Reynolds}

\title{Adaptive time step control for multirate infinitesimal methods \thanks{Submitted to the editors DATE.
\funding{This work was supported by the U.S. Department of Energy, Office of Science, Office of Advanced Scientific Computing Research, Scientific Discovery through Advanced Computing (SciDAC) Program through the FASTMath Institute, under DOE award DE-SC0021354.}}}

\author{
Alex C.~Fish\thanks{Department of Mathematics, Southern Methodist University, Dallas, TX (\email{afish@smu.edu}, \email{reynolds@smu.edu}).}
\and
Daniel R.~Reynolds\footnotemark[2]}

\makeatletter
\newcommand*{\addFileDependency}[1]{
  \typeout{(#1)}
  \@addtofilelist{#1}
  \IfFileExists{#1}{}{\typeout{No file #1.}}
}
\makeatother

\ifpdf
\hypersetup{
  pdftitle={Adaptive time step control for multirate infinitesimal methods},
  pdfauthor={A.~C. Fish and D.~R. Reynolds}
}
\fi

\begin{document}

\maketitle
\begin{abstract}
Multirate methods have been used for decades to temporally evolve initial-value problems in which different components evolve on distinct time scales, and thus use of different step sizes for these components can result in increased computational efficiency.  Generally, such methods select these different step sizes based on experimentation or stability considerations.  For problems that evolve on a single time scale, adaptivity approaches that strive to control local temporal error are widely used to achieve numerical results of a desired accuracy with minimal computational effort, while alleviating the need for manual experimentation with different time step sizes.  However, there is a notable gap in the publication record on the development of adaptive time-step controllers for multirate methods.  In this paper, we extend the single-rate controller work of Gustafsson (1994) to the multirate method setting.  Specifically, we develop controllers based on polynomial approximations to the principal error functions for both the ``fast'' and ``slow'' time scales within \editmade{multirate infinitesimal} (MRI) methods.  We additionally investigate a variety of approaches for estimating the errors arising from each time scale within MRI methods.  We then numerically evaluate the proposed multirate controllers and error estimation strategies on a range of multirate test problems, comparing their performance against an estimated optimal performance.  Through this work, we combine the most performant of these approaches to arrive at a set of multirate adaptive time step controllers that robustly achieve desired solution accuracy with minimal computational effort.
\end{abstract}




\section{Introduction}
\label{sec:intro}
Multirate methods are numerical time integration algorithms used to approximate solutions to ordinary differential equation (ODE) initial-value problems (IVP),
\[y'(t)=f(t,y)=f^{\editmade{\mfs}}(t,y)+f^{\editmade{\mff}}(t,y),\quad y(0)=y_0,\]
in which some portion of the right-hand side function $f(t,y)$ evolves on a slower time scale (and can therefore be evaluated less frequently) than the rest of the function.  Multirate methods typically split $f(t,y)$ additively into two parts, $f^{\editmade{\mfs}}$, the \emph{slow function}, which is evaluated less frequently, and $f^{\editmade{\mff}}$, the \emph{fast function}, which is evaluated more frequently.  For many large-scale applications, multirate methods become particularly attractive when $f^{\editmade{\mfs}}$ has significantly larger computational cost than $f^{\editmade{\mff}}$, and thus a method that minimizes calls to $f^{\editmade{\mfs}}$ may achieve significant gains in computational efficiency for achieving a desired solution accuracy.

A variety of families of multirate methods have been proposed, including MrGARK \cite{gunther_multirate_2016}, MIS \cite{Schlegel2009,Wensch2009}, MRI-GARK \cite{sandu_class_2019}, IMEX-MRI-GARK \cite{chinomona_implicit-explicit_2021}, MERK \cite{luan_new_2020}, MERB \cite{luan_multirate_2021}, and others.  In this work, we focus on so-called \emph{infinitesimal}-type methods, which include all of the above methods except MrGARK.
Generally speaking, an infinitesimal method for evolving a single time step $y_n\approx y(t_n)$ to $y_{n+1}\approx y(t_n+H)$, with \editmade{its embedded solution} $\tilde{y}_{n+1}\approx y(t_n+H)$, proceeds through the following sequence of steps.
\begin{enumerate}
    \item Let: $Y_1=y_n$
    \item For $i=2,...,s$:
      \label{eq:mristep_stages}
      \begin{enumerate}
      \item Solve: $\editmade{v_i}'(\theta) = C_i\, f^{\editmade{\mff}}(\theta,\editmade{v_i}(\theta)) + r_i(\theta)$, for $\theta\in[\theta_{0,i},\theta_{F,i}]$ with $\editmade{v_i}(\theta_{0,i})=v_{0,i}$.\label{eq:mristep_faststage}
      \item Let: $Y_i = \editmade{v_i}(\theta_{F,i})$.
      \end{enumerate}
    \item Solve: $\editmade{\tilde{v}_s}'(\theta) = C_s\, f^{\editmade{\mff}}(\theta,\editmade{\tilde{v}_s}(\theta)) + \tilde{r}_{s}(\theta)$, for $\theta\in[\theta_{0,s},\theta_{F,s}]$ with $\editmade{\tilde{v}_s}(\theta_{0,s})=v_{0,s}$.\label{eq:mristep_fastembedding}
    \item Let: $y_{n+1} = Y_s$ and $\tilde{y}_{n+1} = \editmade{\tilde{v}_s}(\theta_{F,s})$.
\end{enumerate}
In the above, the multirate method is uniquely defined by its choice of leading constant $C_i$, fast stage time intervals $[\theta_{0,i},\theta_{F,i}]$, initial conditions $v_{0,i}$, forcing functions $r_i(\theta)$, and embedding forcing function $\tilde{r}_{s}(\theta)$. These forcing functions are typically constructed using linear combinations of $\left\{f^{\editmade{\mfs}}(\theta_{F,j},Y_j)\right\}$, that serve to propagate information from the slow to the fast time scales. \editmade{The embedded solution $\tilde{y}_{n+1}$ is similar to embedded solutions in standard Runge-Kutta methods. It is a solution with an alternate order of accuracy computed with minimal extra effort once the primary solution $y_{n+1}$ is computed, in this case it is computed by re-solving the last stage with an alternate forcing function, $\tilde{r}_s(\theta)$}.

As seen in step \ref{eq:mristep_faststage} above, computation of each stage in an infinitesimal method requires solving a secondary, \emph{inner}, IVP.
Typically, these inner IVPs are not solved exactly, and instead they are approximated through another numerical method using inner time step size $h \ll H$.  In this work, we consider these inner IVPs to be solved using a one-step method that ``subcycles'' the solution using step sizes $h=\frac{H}{M}$, where $M$ is an integer multirate ratio that describes the difference in dynamical time scales between $f^{\editmade{\mfs}}$ and $f^{\editmade{\mff}}$.

Adaptive time step control is the technique of changing the time step size throughout a solve with the goal of producing a solution which is accurate to a given tolerance of the true solution.
This topic has been thoroughly researched with numerous algorithms developed for single-rate methods, i.e., methods that use a shared step size $H$ for all components of $f(t,y)$ \cite{gustafsson_control-theoretic_1994,soderlind2003,soderlind2006,soderlind2002}.
Savcenco and collaborators developed a component-wise multirate adaptivity strategy designed for use with single-rate methods \cite{savcenco2007}; however, minimal work has been done in developing multirate controllers for use with multirate methods.  Specifically, to our knowledge, the only published work on error-based multiple time step adaptivity for multirate methods was developed by Sarshar and collaborators \cite{sarshar_design_2019}. The two simple adaptivity schemes proposed in that work were designed for MrGARK methods, and were not a primary focus of the paper.

To address the need for a wider ecosystem of algorithms for multirate time step adaptivity, in Section \ref{sec:controllers} we develop a family of simultaneous controllers for both $H$ and $M$, following the techniques of Gustafsson \cite{gustafsson_control-theoretic_1994}, derived using constant and linear approximations of the principal error function.
We additionally introduce two methods that modify our derived controllers so that they more closely emulate the structure of standard single-rate PI \cite[Chapter~IV.2]{kennedy_additive_2003,hairer1996stiff} and PID \cite{kennedy_additive_2003,soderlind2003} controllers.

All of the proposed controllers update both the values of the slow step size $H$ and the multirate ratio $M$ because single-rate adaptivity of just $H$ while keeping $M$ static may lead to an unnecessary amount of computational effort when $H$ is small, or stages with insufficient accuracy when $H$ is large.
Conversely, single-rate adaptivity of only $h=H/M$ while holding $H$ static inhibit the solver flexibility, again either resulting in excessive computational effort or solution error, particularly for problems where the dynamical time scales change throughout the simulation.

Each of our proposed controllers strive to achieve a target overall solution accuracy.  However, errors in \editmade{multirate infinitesimal} approximations can arise from both the choice of $H$ and $M$, where the former determines the approximation error arising from the slow time scale, and the latter determining the error that results from approximating solutions in step \ref{eq:mristep_faststage} of the above algorithm.  We thus propose a set of options for estimating the errors arising from each of these contributions in Section \ref{sec:error_estimation}.


Finally, in Section \ref{sec:results} we identify a variety of test problems and a range of \editmade{multirate infinitesimal} methods to use as a testing suite \editmade{and define metrics with which to measure performance of our controllers}.  We then use this test suite to examine the performance of our proposed algorithms for \editmade{multirate infinitesimal} error estimation, followed by a comparison of the performance of our proposed controllers.


\section{Control-theoretic approaches for multirate adaptivity}
\label{sec:controllers}

The error in \editmade{a multirate infinitesimal} method at the time $t_n$, $\varepsilon_n := \|y(t_n)-y_n\|$, can be bounded by the sum of the \emph{slow error} $\varepsilon_n^{\editmade{\mfs}}$ inherent to the method itself (which assumes that stages in step \ref{eq:mristep_faststage} are solved \emph{exactly}), and the \emph{fast error} $\editmade{\varepsilon_n^{\mff}}$ caused by approximating these stage solutions using some other ``inner'' solver, i.e.,
\[
    \varepsilon_n = \|y(t_n)-y_n\| = \|y(t_n) - y^*_n + y^*_n - y_n\|
    \le \|y(t_n) - y^*_n\| + \|y^*_n - y_n\| = \varepsilon_n^{\editmade{\mfs}} + \varepsilon_n^{\editmade{\mff}}.
\]
Here, we use $y^*_n$ to denote the ``imagined'' solution in which each stage is solved exactly.
We can express the above error contributions over the time step $t_n\to t_{n+1} = t_n+H_n$ as
\begin{equation}
\label{eq:asymptotic_errors}
\begin{aligned}
  \varepsilon^{\editmade{\mfs}}_{n+1} &= \phi^{\editmade{\mfs}}_{n}\, H^P_n \editmade{+ \bigO{H_n^{P+1}}}\\
  \varepsilon^{\editmade{\mff}}_{n+1} &= \phi^{\editmade{\mff}}_{n} \pl\frac{H_n}{M_n}\pr^pH_n \editmade{+ \bigO{\left(\frac{H_n}{M_n}\right)^{p+1}H_n}},
\end{aligned}
\end{equation}
where $P$ and $p$ are the global orders of accuracy for the multirate method embedding and inner method embedding, respectively. We use the embedding orders of accuracy (as opposed to the order of the main methods) since those will be used in Section \ref{sec:error_estimation} for computing our estimates of $\varepsilon^{\editmade{\mfs}}$ and $\varepsilon^{\editmade{\mff}}$.  The above expression for $\varepsilon^{\editmade{\mff}}$ is the global error for a one-step method with step-size $H_n/M_n$ over a time interval of size $H_n$ \cite[Chapter~II.3]{hairer_solving_2009}.  We note the introduction of the subscript ``$n$'' to indicate that these values only apply to the single multirate time step, and that these will be adaptively changed as the simulation proceeds.
The principal error functions $\phi^{\editmade{\mfs}}$ and $\phi^{\editmade{\mff}}$ are functions of $t$ that are unknown in general \editmade{and assumed to be always nonzero}, as they depend on the specific test problem and method used; however, they do not depend on the step size $H_n$ or multirate ratio $M_n$.
\editmade{The notation $\phi_n$ is used in place of $\phi(t_n)$ for brevity. We note that, by convention, $H_n$, $M_n$, $\phi^{\mfs}_n$, and $\phi^{\mff}_n$ give rise to errors $\varepsilon^{\mfs}_{n+1}$ and $\varepsilon^{\mff}_{n+1}$, i.e., at the following time step $t_{n+1}$.}

\editmade{Although stages are typically solved over intervals of size less than $H_n$, the upper bound of the size of the interval is $H_n$, so we use $H_n$ as an approximation for the interval size for the purposes of error estimation. Additionally, the inner method's step size $H_n/M_n$ may be lead to overstepping the end of the interval, in which case the step size for the last step would need to be reduced. We similarly use the upper bound for the inner time step, $H_n/M_n$, as an approximation for the inner time step size for the purposes of error estimation.}

\editmade{In the derivations that follow we treat $M_n$ as a real number, thus when implementing these controllers $M$ must be rounded up to the nearest integer. Dropping the higher order terms from \eqref{eq:asymptotic_errors} and solving them} for $\log(H_n)$ and $\log(M_n)$ gives
\begin{equation}
\begin{aligned}
  \log(H_n) &= \frac{\log(\varepsilon^{\editmade{\mfs}}_{n+1}) - \log(\phi^{\editmade{\mfs}}_{n})}{P} \\
  \log(M_n) &= \frac{(p+1)(\log(\varepsilon^{\editmade{\mfs}}_{n+1})-\log(\phi^{\editmade{\mfs}}_{n}))-P(\log(\varepsilon^{\editmade{\mff}}_{n+1})-\log(\phi^{\editmade{\mff}}_{n})) }{Pp}.
\end{aligned}
\label{eq:solveforloghlogm}
\end{equation}

We wish to choose values of $H_n$ and $M_n$ such that the resulting error $\varepsilon_{n+1}$ is equal to $\tol$, the desired tolerance for the numerical solution to approximate the true solution of the IVP.
Since $\varepsilon_n$ has contributions from both $\varepsilon^{\editmade{\mff}}$ and $\varepsilon^{\editmade{\mfs}}$, we enforce that these \editmade{sources be equal to their time-scale-specific tolerance values}, i.e., we set \editmade{$\varepsilon^{\mfs}_{n+1}=\tol^\mfs$ and $\varepsilon^{\mff}_{n+1}=\tol^\mff$, where $\tol=\tol^\mfs+\tol^\mff$.}
\begin{align*}
  \log(H_n) &= \frac{\log(\editmade{\tol^\mfs}) - \log(\phi^{\editmade{\mfs}}_{n})}{P} \\
  \log(M_n) &= \frac{(p+1)(\editmade{\tol^\mfs})-\log(\phi^{\editmade{\mfs}}_{n}))-P(\log(\editmade{\tol^\mff})-\log(\phi^{\editmade{\mff}}_{n})) }{Pp}.
\end{align*}

Because the principal error function values are generally unknown, we must approximate $\log(\phi^{\editmade{\mfs}}_{n})$ and $\log(\phi^{\editmade{\mff}}_{n})$ with some $\log(\widehat{\phi}^{\editmade{\mfs}}_{n})$ and $\log(\widehat{\phi}^{\editmade{\mff}}_{n})$.
Assuming sufficiently accurate approximations of the principal error functions, we insert these in place of the exact principal error functions above to arrive at our final $\log$-space formulation for the update functions,
\begin{equation}
\begin{aligned}
\log(H_n) &= \frac{\log(\editmade{\tol^\mfs}) - \log(\widehat{\phi}^{\editmade{\mfs}}_{n})}{P} \\
\log(M_n) &= \frac{(p+1)(\log(\editmade{\tol^\mfs})-\log(\widehat{\phi}^{\editmade{\mfs}}_{n}))-P(\log(\editmade{\tol^\mff})-\log(\widehat{\phi}^{\editmade{\mff}}_{n})) }{Pp}.
\label{eq:generalcontroller}
\end{aligned}
\end{equation}

\subsection{Approximations for \texorpdfstring{$\log(\phi(t))$}{log(Phi)}}
\label{sec:phi-approximations}

In \cite{gustafsson_control-theoretic_1994}, Gustafsson developed a single-rate controller by approximating $\log(\phi\editmade{(t)})$ with a piecewise linear function.
We extend this work by approximating $\log(\phi^{\editmade{\mfs}}\editmade{(t)})$ and $\log(\phi^{\editmade{\mff}}\editmade{(t)})$ with piecewise constant functions and reiterate Gustafsson's piecewise linear derivation with further details on motivations throughout.  As with Gustafsson's work, we introduce the shift operator $q$, defined as $q^{a}\log(\phi_n)=\log(\phi_{n+a})$; this is a common operator used in control theory to derive algebraic expressions in terms of a single iteration index \cite{Rota1973}. \editmade{By convention, the shift operator is assumed to be invertible, i.e., $q^{-a}q^a\log(\phi_n)=q^{-a}\log(\phi_{n+a})=\log(\phi_n)$. The identity shift operator is defined as $q^0$ such that $q^{0}\log(\phi_{n})=\log(\phi_n)$.}

\subsubsection{Piecewise Constant \texorpdfstring{$\log(\phi(t))$}{log(Phi)} Approximation}
\label{sec:constant-phi}
If we assume $\log(\phi\editmade{(t)})$ is constant in time, then $\log(\widehat{\phi}_n)=\log(\widehat{\phi}_{n-1})$.
We can convert this to a piecewise constant by adding a corrector based on the true value of $\log(\phi)$, via an ``observer'' in the control \cite[Chapters~11,13]{Dorf2017}.
Introducing the observer with free parameter $k>0$ gives
\begin{align*}
\log(\widehat{\phi}_{n}) &= \log(\widehat{\phi}_{n-1}) + k(\log(\phi_{n-1}) - \log(\widehat{\phi}_{n-1})) \\
&= (1-k)q^{-1}\log(\widehat{\phi}_{n}) + k\log(\phi_{n-1}).
\end{align*}
As in Gustafsson \cite{gustafsson_control-theoretic_1994}, we \editmade{assume that $q+(k-1)q^0$ is invertible to} solve for $\log(\widehat{\phi}_{n})$ to arrive at the piecewise constant approximation,
\begin{align}
    \log(\widehat{\phi}_{n}) &= \frac{kq}{q+\editmade{q^0}(k-1)}\log(\phi_{n-1}).
    \label{eq:piecewiseconstantphi}
\end{align}

\subsubsection{Piecewise Linear \texorpdfstring{$\log(\phi(t))$}{log(Phi)} Approximation}
\label{sec:linear-phi}

In this section we reproduce the derivation of Gustafsson, albeit with additional details that were lacking from \cite{gustafsson_control-theoretic_1994}.  Thus, assuming that $\log(\phi\editmade{(t)})$ is linear (and \editmade{its} time derivative is constant), we have
\begin{align*}
  \log(\widehat{\phi}_{n}) &= \log(\widehat{\phi}_{n-1}) + \nabla\log(\widehat{\phi}_{n-1}), \\
  \nabla\log(\widehat{\phi}_{n}) &= \nabla\log(\widehat{\phi}_{n-1}).
\end{align*}
\editmade{where $\nabla$ is the gradient with respect to the iteration number, indexed by $n$, as in Gustafsson.} We rewrite the system as
\begin{align*}
  \widehat{x}_n &= A \widehat{x}_{n-1} \\
  \log(\widehat{\phi}_n) &= C \widehat{x}_n
\end{align*}
where
\begin{align*}
  \widehat{x}_n &=
  \begin{bmatrix}
    \log(\widehat{\phi}_{n}) \\
    \nabla\log(\widehat{\phi}_{n})
  \end{bmatrix}
  ,\
  \widehat{x}_{n-1} =
  \begin{bmatrix}
    \log(\widehat{\phi}_{n-1}) \\
    \nabla\log(\widehat{\phi}_{n-1})
  \end{bmatrix}
  ,\
  A=
  \begin{bmatrix}
    1 & 1 \\
    0 & 1
  \end{bmatrix},\
  C=
  \begin{bmatrix}
    1 & 0
  \end{bmatrix}.
\end{align*}
Again inserting a corrector based on the true values of $\log(\phi)$ and $\nabla\log(\phi)$, via an observer with free parameter vector $K=[k_1\ k_2]^T$ in the control (here both $k_1,k_2>0$), we have
\begin{align*}
  \widehat{x}_n &= A\widehat{x}_{n-1} + K(\log(\phi_{n-1}) - \log(\widehat{\phi}_{n-1})) \\
  &= (A-KC)\widehat{x}_{n-1} + K\log(\phi_{n-1}).
\end{align*}
Using a backwards difference approximation for $\nabla$, applying the shift operator to convert all approximate $\log(\widehat{\phi})$ terms to the same iteration, summing the equations, and solving for $\log(\widehat{\phi}_n)$ in terms of $\log(\phi_{n-1})$ leads to the piecewise linear approximation for $\log(\phi)$,
\begin{align}
  \log(\widehat{\phi}_n) &= \frac{(k_1+k_2)q^2}{2q^2+q(k_1+k_2-4) +2\editmade{q^0}}\log(\phi_{n-1}),
  \label{eq:piecewiselinearphi}
\end{align}
\editmade{where we again assume that $2q^2+q(k_1+k_2-4) +2q^0$ is invertible.} We note that this approximation will lead to $H$-$M$ controllers with the following structure, which depend on the two most recent values of $H$, but only on the most recent error term:
\begin{align*}
    H_{n+1}&=H_n^{\gamma_1}H_{n-1}^{\gamma_2}\left(\frac{\tol}{\varepsilon^{\editmade{\mfs}}_{n+1}}\right)^{\alpha}
\end{align*}
This short error history was expected for our earlier piecewise constant derivation, but, because the controller now depends on a longer history of $H$, it is typical to convert this to depend on an equal history of error terms.
We follow Gustafsson's derivation by enforcing the additional condition
\begin{align}
  k_1\log(\widehat{\phi}_{n-1}) = k_1\log(\phi_{n-1})
  \label{eq:piecewiselinearcorrection}
\end{align}
that enforces continuity in the piecewise approximation.  Subtracting this from \eqref{eq:piecewiselinearphi}, we may arrive at the modified approximation for $\log(\phi_n)$, \editmade{assuming $2q^2+q(k_1+k_2-4) +q^0(2-k_1)$ is invertible,}
\begin{equation}
  \log(\widehat{\phi}_n) = \frac{(k_1+k_2)q^2-k_1q}{2q^2+q(k_1+k_2-4) +\editmade{q^0}(2-k_1)}\log(\phi_{n-1}).
  \label{eq:piecewiselinearphicorrected}
\end{equation}

\subsubsection{Expression for
\texorpdfstring{$\log(\phi_{n-1})$}{log(Phi)}}
\label{sec:exactlogphiprevious}

Both the piecewise constant and piecewise quadratic approaches above resulted in approximations \eqref{eq:piecewiseconstantphi} and \eqref{eq:piecewiselinearphicorrected} that depend on the true values of the principal error function at the previous step, $\log(\phi_{n-1})$.
Although these principal error function values are generally unknown, they may be approximated using estimates of the error at the end of the previous step.  Depending on whether $\log(\phi^{\editmade{\mfs}})$ or $\log(\phi^{\editmade{\mff}})$ is being approximated, we may estimate the principal error function at the previous time step by rearranging \eqref{eq:solveforloghlogm},
\begin{equation}
\begin{aligned}
\log(\phi^{\editmade{\mfs}}_{n-1}) &= \log(\varepsilon^{\editmade{\mfs}}_n) - Pq^{-1}\log(H_n) \\
\log(\phi^{\editmade{\mff}}_{n-1}) &= \log(\varepsilon^{\editmade{\mff}}_n) - (p+1)q^{-1}\log(H_n) + q^{-1}p\log(M_n),
\label{eq:truelogphin}
\end{aligned}
\end{equation}

\subsection{H-M Controllers}
\label{sec:hm-controllers}

At this point we have all of the requisite parts to finalize controller formulations for $H_n$ and $M_n$.  To this end, we insert the expressions \eqref{eq:truelogphin} into the approximations \eqref{eq:piecewiseconstantphi} or \eqref{eq:piecewiselinearphicorrected}.
We then insert the resulting formulas into the $\log$-space formulation of the update functions \eqref{eq:generalcontroller}, apply the shift operator, and solve the resulting equations for $\log(H_{n+1})$ and $\log(M_{n+1})$ in terms of known or estimated values of $\log(H_n)$, $\log(M_n)$, $\log(\varepsilon^{\editmade{\mfs}}_n)$, and $\log(\varepsilon^{\editmade{\mff}}_n)$.
We then finalize each controller by taking the exponential of these log-space update functions.

For simplicity of notation, we introduce the ``oversolve'' factors,
$\eta^{\editmade{\mfs}}_n := \left(\editmade{\tol^\mfs}\right)/\varepsilon^{\editmade{\mfs}}_n$ and $\eta^{\editmade{\mff}}_n := \left(\editmade{\tol^\mff}\right)/\varepsilon^{\editmade{\mff}}_n$.  The values of these factors reflect the deviation between the achieved errors and their target, where a value of one indicates an optimal control.

There are four potential combinations of our two piecewise approximations for each of $\log(\phi^{\editmade{\mfs}}_{n})$ and $\log(\phi^{\editmade{\mff}}_{n})$. While using differing order approximations is valid in theory, we have not found this useful in practice.
In particular, due to the coupled nature of the general controllers \eqref{eq:generalcontroller}, the \editmade{extent of the error history for both fast \emph{and} slow errors in the $M_n$ update function will correspond to the higher-order of the $\log(\phi^{\mfs}_{n})$ and $\log(\phi^{\mff}_{n})$ approximations,} and thus it makes little sense to use a lower-order approximation for $\log(\phi^{\editmade{\mff}}_{n})$ than $\log(\phi^{\editmade{\mfs}}_{n})$.  Additionally, we expect the fast time scale (associated with $M$) to vary more rapidly than the slow time scale for multirate applications, so we also we rule out cases where the approximation order for $\log(\phi^{\editmade{\mff}}_{n})$ exceeds \editmade{that of $\log(\phi^{\mfs}_{n})$}.

\subsubsection{Constant-Constant Controller}
\label{sec:cc-controller}

When using a piecewise constant approximation for both $\log(\phi^{\editmade{\mfs}})$ and $\log(\phi^{\editmade{\mff}})$, we follow the procedure outlined above to arrive at the following equations for $\log(H_n)$ and $\log(M_n)$,
\begin{align*}
    \log(H_n) &= \frac{k_1q}{P(q-1)}\log\left(\eta_{n}^{\editmade{\mfs}}\right), \\
    \log(M_n) &= \frac{(p+1)k_1q}{Pp(q-1)} \log\pl\eta_{n}^{\editmade{\mfs}}\pr - \frac{k_2q}{p(q-1)}\log\pl\eta_{n}^{\editmade{\mff}}\pr.
\end{align*}
These give the controller pair
\begin{equation}
\label{eq:constant_constant_controller}
\begin{aligned}
    H_{n+1} &= H_n\pl\eta^{\editmade{\mfs}}_{n+1}\pr^\alpha, \\
    M_{n+1} &= M_n\pl\eta^{\editmade{\mfs}}_{n+1}\pr^{\beta_1}\pl\eta^{\editmade{\mff}}_{n+1}\pr^{\beta_2},
\end{aligned}
\end{equation}
where
\[
   \alpha = \dfrac{k_1}{P}, \quad \beta_1 = \dfrac{(p+1)k_1}{Pp}, \quad\text{and}\quad \beta_2=\dfrac{-k_2}{p}.
\]

\subsubsection{Linear-Linear Controller}
\label{sec:ll-controller}

Similarly, using a piecewise linear approximation for both $\log(\phi^{\editmade{\mfs}})$ and $\log(\phi^{\editmade{\mff}})$, and solving for $\log(H_n)$ and $\log(M_n)$, we have
\begin{align*}
    \log(H_n) &= \frac{( k_{11}+k_{12}) q^2-k_{11}q}{2P(q^2-2q+1)}\log\pl\eta^{\editmade{\mfs}}_n\pr \\
    \log(M_n) &= \frac{\pl p+1\pr\pl\pl k_{11}+k_{12}\pr q^2-k_{11}q\pr}{2Pp(q^2-2q+1)}\log\pl\eta^{\editmade{\mfs}}_n\pr \\
    &\quad-  \frac{(k_{21}+k_{22})q^2 - k_{21}q}{2p(q^2-2q+1)}\log\pl\eta^{\editmade{\mff}}_n\pr.
\end{align*}
This results in the controller pair
\begin{equation}
\label{eq:linear_linear_controller}
\begin{aligned}
    H_{n+1} &= H_n^2H_{n-1}^{-1}
    \pl\eta^{\editmade{\mfs}}_{n+1}\pr^{\alpha_1}
    \pl\eta^{\editmade{\mfs}}_{n}\pr^{\alpha_2}, \\
    M_{n+1} &= M_n^2M_{n-1}^{-1}
    \pl\eta^{\editmade{\mfs}}_{n+1}\pr^{\beta_{11}}
    \pl\eta^{\editmade{\mfs}}_{n}\pr^{\beta_{12}}
    \pl\eta^{\editmade{\mff}}_{n+1}\pr^{\beta_{21}}
    \pl\eta^{\editmade{\mff}}_{n}\pr^{\beta_{22}},
\end{aligned}
\end{equation}
where
\begin{equation}
\label{eq:linear_linear_parameters}
\begin{aligned}
    \alpha_1 &= \frac{k_{11}+k_{12}}{2P}, & \alpha_2 &= \frac{-k_{11}}{2P} \\
    \beta_{11} &= \frac{(p+1)(k_{11}+k_{12})}{2Pp}, & \beta_{12} &= \frac{-(p+1)k_{11}}{2Pp}, \\
    \beta_{21} &= \frac{-(k_{21}+k_{22})}{2p}, & \beta_{22} &= \frac{k_{21}}{2p}.
\end{aligned}
\end{equation}

\subsection{A Note on Higher Order \texorpdfstring{$\log(\phi(t))$}{log(Phi)} Approximations}

Following the preceding approaches, it is possible to form controllers based on higher order polynomial approximations to $\log(\phi\editmade{(t)})$.
While we have investigated these, we did not find such multirate $H$-$M$ controllers to have robust performance.  For example, a ``Quadratic-Quadratic'' controller that uses piecewise quadratic approximation for both $\log(\phi^{\editmade{\mfs}}_n)$ and $\log(\phi^{\editmade{\mff}}_n)$ depends on an extensive history of both $H$ or $M$ terms, and by design attempts to limit changes in $H$ and $M$ to ensure that these vary smoothly throughout a simulation.  In our experiments, we found that this over-accentuation of smoothness rendered the controller to lack sufficient flexibility to grow or shrink $H$ or $M$ at a rate at which some multirate problems demanded.
The Quadratic-Quadratic controller thus exhibited a ``bimodal'' behavior, in that for some problems it failed quickly, whereas for others it performed excellently, but we were unable to tune it so that it would provide robustness across our test problem suite.
Since our goal in this paper is to construct controllers with wide applicability to many multirate applications, we omit the Quadratic-Quadratic controller, or any other controllers derived from a higher order $\log(\phi_n)$ approximations, from this paper.

However, to address the over-smoothness issue we introduce two additional controllers, one that results from a simple modification of the Linear-Linear controller above, and the other that extends this idea to include additional error estimates.

\subsection{PIMR Controller}

The ``PI'' controller, popular in single-rate temporal adaptivity, has a similar qualitative structure to the Linear-Linear controller \eqref{eq:linear_linear_controller}, albeit with a reduced dependence on the history of $H_n$ as in \cite{kennedy_additive_2003} and  \cite[Chapter~IV.2]{hairer1996stiff},
\begin{equation}
    \label{eq:PI_controller}
    H_{n+1} = H_n \left(\frac{\text{tol}}{\varepsilon_{n+1}}\right)^{\alpha_1} \left(\frac{\text{tol}}{\varepsilon_{n}}\right)^{\alpha_2}.
\end{equation}
We thus introduce a ``PIMR'' controller by eliminating one factor of both $H_n/H_{n-1}$ and $M_n/M_{n-1}$ from \eqref{eq:linear_linear_controller}:
\begin{equation}
\label{eq:pimr_controller}
\begin{aligned}
    H_{n+1} &= H_n
    \pl\eta^{\editmade{\mfs}}_{n+1}\pr^{\alpha_1}
    \pl\eta^{\editmade{\mfs}}_{n}\pr^{\alpha_2}, \\
    M_{n+1} &= M_n
    \pl\eta^{\editmade{\mfs}}_{n+1}\pr^{\beta_{11}}
    \pl\eta^{\editmade{\mfs}}_{n}\pr^{\beta_{12}}
    \pl\eta^{\editmade{\mff}}_{n+1}\pr^{\beta_{21}}
    \pl\eta^{\editmade{\mff}}_{n}\pr^{\beta_{22}},
\end{aligned}
\end{equation}
where $\alpha_1$, $\alpha_2$, $\beta_{11}$, $\beta_{12}$, $\beta_{21}$ and $\beta_{22}$ are defined as in \eqref{eq:linear_linear_parameters}.

\subsection{PIDMR Controller}

The higher-order ``PID'' controller from single-rate temporal adaptivity has a similar qualitative structure to the PI controller, but includes dependence on one additional error term \cite{kennedy_additive_2003,soderlind2003},
\begin{equation}
    \label{eq:PID_controller}
    H_{n+1} = H_n \left(\frac{\text{tol}}{\varepsilon_{n+1}}\right)^{\alpha_1} \left(\frac{\text{tol}}{\varepsilon_{n}}\right)^{\alpha_2} \left(\frac{\text{tol}}{\varepsilon_{n-1}}\right)^{\alpha_3}.
\end{equation}
The parameters in PID controllers typically alternate signs, with $\alpha_1>0$, $\alpha_2<0$ and $\alpha_3>0$, where each takes the form of a constant divided by the asymptotic order of accuracy for the single-rate method.

We thus introduce an extension of the PIMR controller to emulate this ``PID'' structure, resulting in the ``PIDMR'' controller.
We strive to make this extension as natural as possible -- we thus introduce new free parameters $k_{13}$ and $k_{23}$, adjust the coefficient in the denominator of the exponents from 2 to 3, and assume that the exponents alternate signs.  We thus define the PIDMR controller as
\begin{equation}
\label{eq:pidmr_controller}
\begin{aligned}
    H_{n+1} &= H_n
    \pl\eta^{\editmade{\mfs}}_{n+1}\pr^{\alpha_1}
    \pl\eta^{\editmade{\mfs}}_{n}\pr^{\alpha_2}
    \pl\eta^{\editmade{\mfs}}_{n-1}\pr^{\alpha_3}, \\
    M_{n+1} &= M_n
    \pl\eta^{\editmade{\mfs}}_{n+1}\pr^{\beta_{11}}
    \pl\eta^{\editmade{\mfs}}_{n}\pr^{\beta_{12}}
    \pl\eta^{\editmade{\mfs}}_{n-1}\pr^{\beta_{13}}
    \pl\eta^{\editmade{\mff}}_{n+1}\pr^{\beta_{21}}
    \pl\eta^{\editmade{\mff}}_{n}\pr^{\beta_{22}}
    \pl\eta^{\editmade{\mff}}_{n-1}\pr^{\beta_{23}},
\end{aligned}
\end{equation}
where
\begin{align*}
    \alpha_1 &= \frac{k_{11}+k_{12}+k_{13}}{3P}, & \alpha_2 &= \frac{-(k_{11}+k_{12})}{3P}, & \alpha_3 &= \frac{k_{11}}{3P},\\
    \beta_{11} &= \frac{(p+1)(k_{11}+k_{12}+k_{13})}{3Pp}, & \beta_{12} &= \frac{-(p+1)(k_{11}+k_{12})}{3Pp}, & \beta_{13} &= \frac{(p+1)k_{11}}{3Pp}, \\
    \beta_{21} &= \frac{-(k_{21}+k_{22}+k_{23})}{3p}, & \beta_{22} &= \frac{k_{21}+k_{22}}{3p}, & \beta_{23} &= \frac{-k_{21}}{3p}.
\end{align*}



\section{\editmade{Multirate infinitesimal} method error estimation}
\label{sec:error_estimation}

In infinitesimal methods, the embedded solution is only computed in step \ref{eq:mristep_fastembedding}, i.e., as with standard Runge--Kutta methods the embedding \emph{reuses} the results of each internal stage computation (step \ref{eq:mristep_stages}).
Thus to measure $\editmade{\varepsilon_{\editmade{n+1}}^{\mfs}}$ we will use a standard last-stage embedding of the infinitesimal method, e.g., $\varepsilon_{\editmade{n+1}}^{\editmade{\mfs}} \approx \|y_{\editmade{n+1}}-\tilde{y}_{\editmade{n+1}}\|$.  However, estimation of $\varepsilon_{\editmade{n+1}}^{\editmade{\mff}}$ is less obvious.

Within each of the fast IVP solves in steps \ref{eq:mristep_faststage} and \ref{eq:mristep_fastembedding}, one must employ a separate IVP solver.  For these, we assume that an embedded Runge--Kutta method is used which provides two distinct solutions of differing orders of accuracy for the fast IVP solution.  There are a number of ways to utilize the embedded fast Runge--Kutta method to estimate $\editmade{\varepsilon_{\editmade{n+1}}^{\mff}}$.  We outline five potential approaches with a range of computational costs here, and will compare their performance in Section \ref{sec:results}.

\subsection{Full-Step (FS) strategy}
\begin{figure}[htb]
    \centering
    \includegraphics[width=0.8\textwidth]{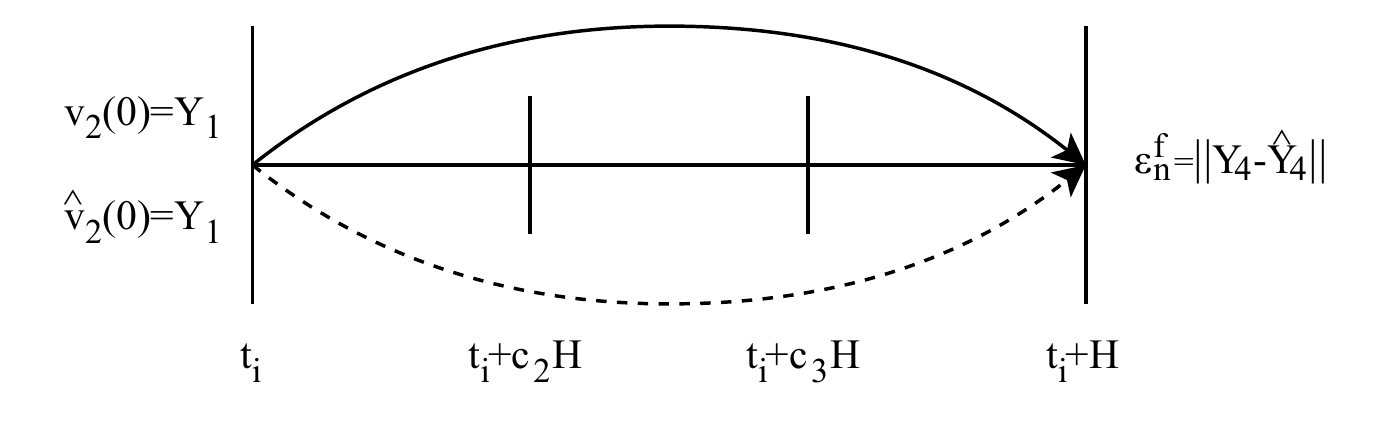}
    \caption{\editmade{\emph{Full-Step} strategy for fast error estimation. A full step is solved with both the primary fast method (solid line) and the embedded fast method (dashed line), which start from the same initial condition $Y_1$ at the first stage solve.}}
    \label{fig:fs-diagram}
\end{figure}

Our first approach for estimating the fast time scale error runs the entire infinitesimal step twice, once with a primary fast method for all fast IVP solves (producing the solution $y_{\editmade{n+1}}$) and once with the fast method's embedding (to produce the solution $\widehat{y}_{\editmade{n+1}}$).  We use the difference of these to estimate the fast error, $\varepsilon^{\editmade{\mff}}_{\editmade{n+1}}\approx\|y_{\editmade{n+1}}-\widehat{y}_{\editmade{n+1}}\|$:
\begin{enumerate}
    \item Let: $Y_1=y_n$
    \item For $i=2,...,s$:
      \begin{enumerate}
      \item Solve with primary fast method: $\editmade{v_i}'(\theta) = C_i\, f^{\editmade{\mff}}(\theta,\editmade{v_i}(\theta)) + r_i(\theta)$, for $\theta\in[\theta_{0,i},\theta_{F,i}]$ with $\editmade{v_i}(\theta_{0,i})=v_{0,i}$.
      \item Solve with embedded fast method: $\editmade{\widehat{v}_i}'(\theta) = C_i\, f^{\editmade{\mff}}(\theta,\editmade{\widehat{v}_i}(\theta)) + \widehat{r}_i(\theta)$, for $\theta\in[\theta_{0,i},\theta_{F,i}]$ with $\editmade{\widehat{v}_i}(\theta_{0,i})=\widehat{v}_{0,i}$.
      \label{eq:fullstep_embedding}
      \item Let: $Y_i = \editmade{v_i}(\theta_{F,i})$ and $\widehat{Y_i} = \editmade{\widehat{v}_i}(\theta_{F,i})$.
      \end{enumerate}
    \item Let: $y_{n+1} = Y_s$ and $\widehat{y}_{n+1} = \widehat{Y}_s$.
    \item Let: $\varepsilon^{\editmade{\mff}}_{\editmade{n+1}} = \|y_{n+1}-\widehat{y}_{n+1}\|$.
\end{enumerate}
While this is perhaps the most accurate method for estimating $\varepsilon^{\editmade{\mff}}$, it requires two separate computations of the full infinitesimal step, including two sets of evaluations of the slow right-hand side function, $f^{\editmade{\mfs}}$. \editmade{See Figure \ref{fig:fs-diagram} for a diagram of this strategy.}

\subsection{Stage-Aggregate (SA) strategy}
\begin{figure}[htb]
    \centering
    \includegraphics[width=0.95\textwidth]{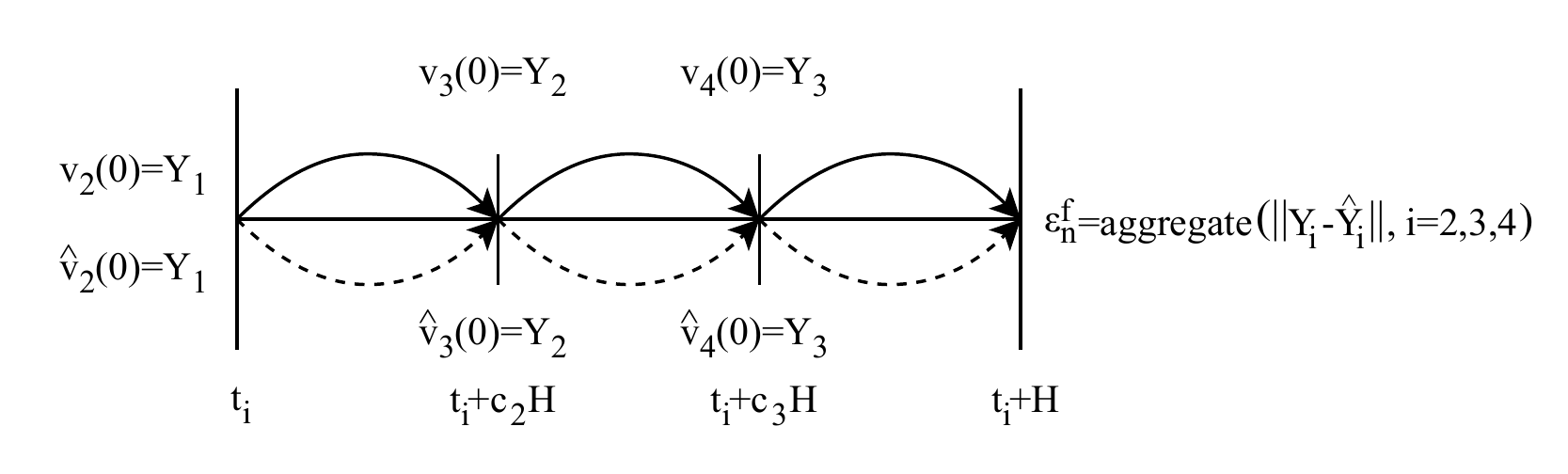}
    \caption{\editmade{\emph{Stage-Aggregate} strategy for fast error estimation. Each stage is solved with both the primary fast method (solid line) and the embedded fast method (dashed line). Both the primary and embedded fast methods use the result of the primary fast method as their initial condition for the next stage solve. The error is computed with an aggregating function on the stage errors.}}
    \label{fig:sa-diagram}
\end{figure}

As above, we consider two separate fast IVP solves for each stage, once with the primary fast method and once with its embedding.  However, here both fast IVP solves use the results from the \textit{primary} fast method for their initial conditions and the forcing terms $r_i(\theta)$.  Here, we estimate the fast error by aggregating the norms of the difference in computed stages, denoted as
\[
  \varepsilon^{\editmade{\mff}}_{\editmade{n+1}}\approx\textnormal{aggregate}(\|Y_i-\widehat{Y}_i\|,\; i=2,\ldots,s),
\]
where $Y_i$ are the stages computed using the primary fast method, and $\widehat{Y}_i$ are the stages computed using the fast method's embedding.  We consider aggregating by mean and max, and refer to these as \emph{Stage-Aggregate-mean} (SA-mean) and \emph{Stage-Aggregate-max} (SA-max), respectively:
\begin{enumerate}
    \item Let: $Y_1=y_n$
    \item For $i=2,...,s$:
      \begin{enumerate}
      \item Solve once using primary fast method and once using embedded fast method: $\editmade{v_i}'(\theta) = C_i\, f^{\editmade{\mff}}(\theta,\editmade{v_i}(\theta)) + r_i(\theta)$, for $\theta\in[\theta_{0,i},\theta_{F,i}]$ with $\editmade{v_i}(\theta_{0,i})=v_{0,i}$.
      \item Let: $Y_i = \editmade{v_i}(\theta_{F,i})$ when using primary fast method.
      \item Let: $\widehat{Y}_i = \editmade{v_i}(\theta_{F,i})$ when using embedded fast method.
      \end{enumerate}
    \item Let: $y_{n+1}=Y_s$.
    \item Let: $\varepsilon^{\editmade{\mff}}_{\editmade{n+1}} = \textnormal{aggregate}(\|Y_i-\widehat{Y}_i\|,\; i=2,\ldots,s)$.
\end{enumerate}
Because this estimation strategy solves the \emph{same} IVP at each stage using two separate methods, the number of $f^{\editmade{\mfs}}$ evaluations is cut in half compared to the \emph{Full-Step} strategy. \editmade{See Figure \ref{fig:sa-diagram} for a diagram of this strategy.}

\subsection{Local-Accumulation-Stage-Aggregate (LASA) strategy}
\begin{figure}[htb]
    \centering
    \includegraphics[width=1.0\textwidth]{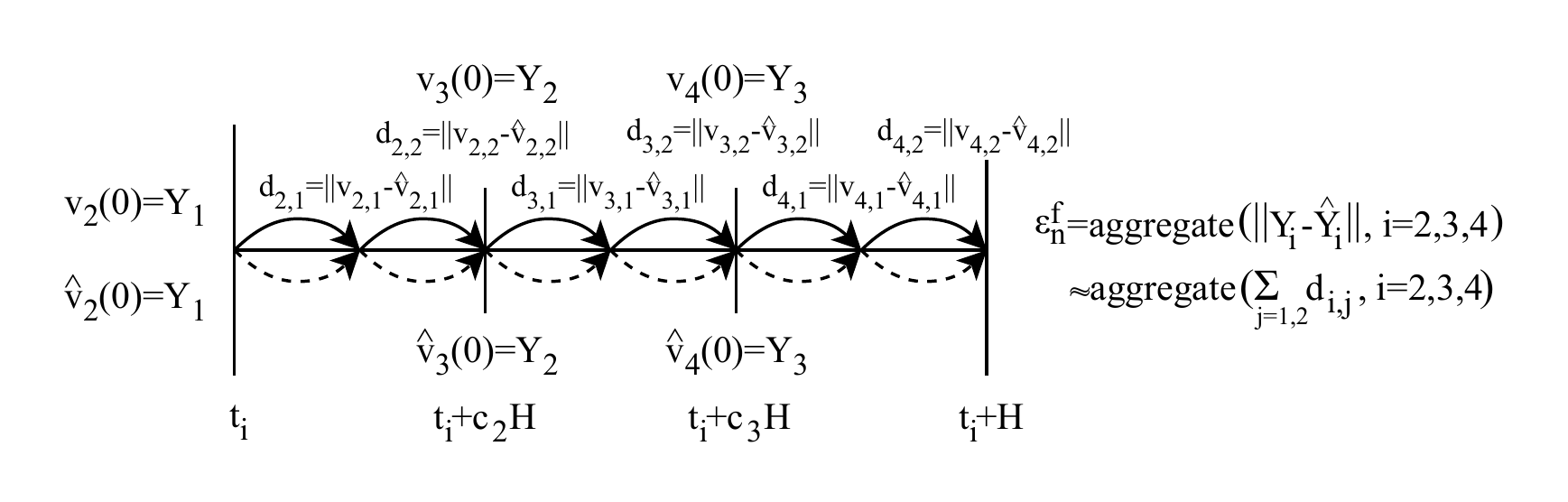}
    \caption{\editmade{\emph{Local-Accumulation-Stage-Aggregate} strategy for fast error estimation. Each step of each stage is solved with both the primary fast method (solid line) and the embedded fast method (dashed line), requiring no extra function evaluations. Each subsequent step of the fast method uses the result from the primary fast method as its initial condition. The norm of the difference between the fast step solutions is summed for each multirate stage and used as an approximation of the stage error, and the overall error is computed by aggregating the stage error approximations.}}
    \label{fig:lasa-diagram}
\end{figure}

Our final $\varepsilon^{\editmade{\mff}}_{\editmade{n+1}}$ measurement strategy is designed to minimize the overall algorithmic cost.  As with the \emph{Stage-Aggregate} approaches above, we solve a single set of fast IVPs, however instead of using the embedded method solution as a separate fast solver, we evolve only the primary fast method and use its embedding to estimate the temporal error within each fast sub-step.  To estimate the overall error within each slow stage, we sum the fast sub-step errors, i.e., for the slow stage $i$ and fast step $j$, we define $d_{j,i}=\|v_{j,i}-\widehat{v}_{j,i}\|$ as the norm of the difference between the primary fast method and its embedding.  Then as with the \emph{Stage-Aggregate} strategy, we aggregate the stage errors to estimate
\begin{align*}
    \varepsilon^{\editmade{\mff}}_{\editmade{n+1}} &\approx\textnormal{aggregate}(\|Y_i-\widehat{Y}_i\|,\;i=2,\ldots,s) \\
    &\approx\textnormal{aggregate}\left(\sum_{j=1}^M d_{j,i},\;i=2,\ldots,s\right).
\end{align*}
We again consider aggregating by mean and max, and we refer to the corresponding estimation strategies as \emph{Local-Accumulation-Stage-Aggregate-mean} (LASA-mean) and \emph{Local-Accumulation-Stage-Aggregate-max} (LASA-max), respectively:
\begin{enumerate}
    \item Let: $Y_1=y_n$
    \item For $i=2,...,s$:
      \begin{enumerate}
      \item Solve : $\editmade{v_i}'(\theta) = C_i\, f^{\editmade{\mff}}(\theta,\editmade{v_i}(\theta)) + r_i(\theta)$, for $\theta\in[\theta_{0,i},\theta_{F,i}]$ with $\editmade{v_i}(\theta_{0,i})=v_{0,i}$.
      \begin{enumerate}
          \item Let: $v_{j,i}$ be the step solution using the primary fast method at sub-step $j$, $j=1,...,M$.
          \item Let: $\widehat{v}_{j,i}$ be the step solution using the embedded fast method at sub-step $j$, $j=1,...,M$.
          \item Let: $d_{j,i}=\|v_{j,i}-\widehat{v}_{j,i}\|$.
          \item Use $v_{j,i}$ as the initial condition for the next step.
      \end{enumerate}
      \item Let: $Y_i = \editmade{v_i}(\theta_{F,i})$
      \item Let: $\|Y_i-\widehat{Y}_i\| \approx \sum_{j=1}^M d_{j,i}$.
      \end{enumerate}
    \item Let: $y_{n+1}=Y_s$.
    \item Let: $\varepsilon^{\editmade{\mff}}_{\editmade{n+1}} = \textnormal{aggregate}(\|Y_i-\widehat{Y}_i\|,\; i=2,\ldots,s)$.
\end{enumerate}
Here, since all fast IVPs use the same set of $f^{\editmade{\mfs}}$ evaluations, then this approach uses half as many slow right-hand side evaluations as the \emph{Full-Step} strategy.  Additionally, since each fast IVP is solved only once, then it uses approximately half the number of $f^{\editmade{\mff}}$ evaluations in comparison to both the \emph{Full-Step} and \emph{Stage-Aggregate} strategies. \editmade{See Figure \ref{fig:lasa-diagram} for a diagram of this strategy.}

In Section \ref{sec:esf-strategies}, we compare the performance of these five strategies, based on both how close their solutions come to the desired tolerance, and on the relative computational cost of each strategy.

\section{Numerical results}
\label{sec:results}

We assess performance for our five measurement strategies and four multirate controllers.  We base these assessments both with respect to how well their resulting solutions match the desired tolerance, and also how close their computational costs were to \editmade{a set of estimated optimal} costs\editmade{, which are discussed in Section \ref{sec:assessing-adaptive-performance} and Appendix \ref{appendix:optimal-performance}}. All of the codes used for these computational results are available in the public GitHub repository:\\ \hyperlink{https://github.com/fishac/AdaptiveHMControllerPaper}{https://github.com/fishac/AdaptiveHMControllerPaper}.


\subsection{Test problems}
\label{sec:results-test_problems}

We assessed the performance of both our fast error estimation strategies and multirate controllers using a test suite comprising seven multirate test problems.  For each problem, we compute error using analytical solutions when available; otherwise we use MATLAB's \emph{ode15s} with tight tolerances $AbsTol=10^{-14}$ and $RelTol=2.5\times 10^{-14}$ to compute reference solutions at \editmade{ten evenly spaced points in the problem's time interval}.

\subsubsection{Bicoupling}
\label{problem:bicoupling}

The Bicoupling problem is a nonlinear and nonautonomous multirate test problem proposed in \cite{luan_multirate_2021},
\begin{align*}
    \begin{bmatrix}
      u' \\
      v' \\
      w'
    \end{bmatrix}
    =
    \begin{bmatrix}
      \gamma v - w - pt \\
      -\gamma u \\
      -lw-lpt - p \pl u-\dfrac{aw}{al+b\gamma} - \dfrac{apt}{al+b\gamma}\pr^2 - p \pl v-\dfrac{bw}{al+b\gamma} - \dfrac{bpt}{al+b\gamma}\pr^2
    \end{bmatrix},
\end{align*}
for $t\in[0,1]$, with initial conditions $u(0) = 1+a$, $v(0) = b$, $w(0) = al+b\gamma$ and parameters $a=1$, $b=20$, $\gamma=100$, $l=5$, and $p=0.01$.  This IVP has true solution
\begin{align*}
    u(t) = \cos(\gamma t) + ae^{-lt}, \quad
    v(t) = -\sin(\gamma t) + be^{-lt}, \quad
    w(t) = (al+b\gamma)e^{-lt}-pt.
\end{align*}

We apply the same multirate splitting of the right-hand side function as \cite{luan_multirate_2021}, that used
$f^{\editmade{\mfs}} = \begin{bmatrix} \gamma v & -\gamma u & 0 \end{bmatrix}^T$ and $f^{\editmade{\mff}} = f - f^{\editmade{\mfs}}$.
%

\subsubsection{Stiff Brusselator ODE}
\label{problem:brusselator}

The Brusselator is an oscillating chemical reaction problem which is widely used to test multirate, implicit, and mixed implicit-explicit methods.
We define the problem with the same parameters as in \cite{luan_new_2020},
\begin{align*}
    \begin{bmatrix}
      u' \\
      v' \\
      w'
    \end{bmatrix}
    &=
    \begin{bmatrix}
      a - (w+1)u + u^2v \\
      uw-u^2v \\
      \dfrac{b-w}{\epsilon} - uw
    \end{bmatrix},\qquad
    \begin{bmatrix}
      u(0)\\
      v(0) \\
      w(0)
    \end{bmatrix}
    =
    \begin{bmatrix}
      1.2 \\
      3.1 \\
      3
    \end{bmatrix}
\end{align*}
for $t\in[0,2]$, and using the parameters $a=1$, $b=3.5$, and $\varepsilon=0.01$.
As we are unaware of an analytical solution to this IVP, we compute reference solutions as described above.

We apply the same multirate splitting of the right-hand side function as in \cite{luan_new_2020},
\begin{align*}
    f^{\editmade{\mfs}} &=
    \begin{bmatrix}
      a + (w+1)u + u^2v \\
      uw-u^2v \\
      \dfrac{b}{\epsilon} - uw
    \end{bmatrix},\qquad
    f^{\editmade{\mff}} =
    \begin{bmatrix}
      0 \\
      0 \\
      -\dfrac{w}{\epsilon}
    \end{bmatrix}.
\end{align*}

\subsubsection{Kaps}
\label{problem:kaps}

The Kaps problem is an autonomous nonlinear problem with analytical solution that has been frequently used to test Runge--Kutta methods, presented in \cite{Skvortsov2003},
\begin{align*}
    \begin{bmatrix}
      u' \\
      v'
    \end{bmatrix}
    &=
    \begin{bmatrix}
      -(\mu+2)u+\mu v^2 \\
      -v^2+u-v
    \end{bmatrix},\qquad
    \begin{bmatrix}
      u(0)\\
      v(0)
    \end{bmatrix}
    =
    \begin{bmatrix}
      1 \\
      1
    \end{bmatrix}
\end{align*}
for $t\in[0,2]$, where we use the stiffness parameter $\mu=100$.  This IVP has true solution
\begin{align*}
    u(t) = e^{-2t}, \qquad v(t) = e^{-t},
\end{align*}
and we split the right-hand side function into slow and fast components as
\begin{align*}
    f^{\editmade{\mfs}} =
    \begin{bmatrix}
      0 \\
      -v^2+u-v
    \end{bmatrix},\qquad
    f^{\editmade{\mff}} =
    \begin{bmatrix}
      -(\mu+2)u+\mu v^2 \\
      0
    \end{bmatrix}.
\end{align*}

\subsubsection{KPR}
\label{problem:kpr}

The KPR problem is a nonlinear IVP system with analytical solution, with variations that have been widely applied to test multirate algorithms.  We use the same formulation as in \cite{chinomona_implicit-explicit_2021},
\begin{align*}
    \begin{bmatrix}
      u' \\
      v'
    \end{bmatrix}
    &=
    \Lambda
    \begin{bmatrix}
      \frac{-3+u^2-\cos(\beta t)}{2u} \\
      \frac{-2+v^2-\cos(t)}{2v}
    \end{bmatrix}
    -
    \begin{bmatrix}
      \frac{\beta\sin(\beta t)}{2u} \\
      \frac{\sin(t)}{2v}
    \end{bmatrix},\qquad
    \begin{bmatrix}
      u(0)\\
      v(0)
    \end{bmatrix}
    =
    \begin{bmatrix}
      2 \\
      \sqrt{3}
    \end{bmatrix}
    \\
    \Lambda &=
    \begin{bmatrix}
      \lambda^{\editmade{\mff}} & \frac{1-\varepsilon}{\alpha}(\lambda^{\editmade{\mff}}-\lambda^{\editmade{\mfs}}) \\
      -\alpha\varepsilon(\lambda^{\editmade{\mff}}-\lambda^{\editmade{\mfs}}) & \lambda^{\editmade{\mfs}}
    \end{bmatrix},
\end{align*}
for $t\in[0, 5\pi/2]$, and with the parameters $\lambda^{\editmade{\mfs}}=-1$, $\lambda^{\editmade{\mff}}=-10$, $\alpha=1$, $\beta=20$, and $\epsilon=0.1$.
This IVP has true solution
\begin{align*}
    u(t) = \sqrt{3+\cos(\beta t)}, \qquad v(t) = \sqrt{2+\cos(t)},
\end{align*}
and we split the right-hand side function component-wise as in \cite{chinomona_implicit-explicit_2021,sandu_class_2019},
\begin{align*}
    f^{\editmade{\mfs}} &=
    \begin{bmatrix}
      0 & 0 \\
      0 & 1
    \end{bmatrix}
    \Lambda
    \begin{bmatrix}
      \frac{-3+u^2-\cos(\beta t)}{2u} \\
      \frac{-2+v^2-\cos(t)}{2v}
    \end{bmatrix}
    -
    \begin{bmatrix}
      0 \\
      \frac{\sin(t)}{2v}
    \end{bmatrix},\\
    f^{\editmade{\mff}} &=
    \begin{bmatrix}
      1 & 0 \\
      0 & 0
    \end{bmatrix}
    \Lambda
    \begin{bmatrix}
      \frac{-3+u^2-\cos(\beta t)}{2u} \\
      \frac{-2+v^2-\cos(t)}{2v}
    \end{bmatrix}
    -
    \begin{bmatrix}
      \frac{\beta\sin(\beta t)}{2u} \\
      0
    \end{bmatrix}.
\end{align*}

\subsubsection{Forced Van der Pol}
\label{problem:forcedvanderpol}

We consider a forced version of the widely-used Van der Pol oscillator test problem, defined in \cite{DBLP:journals/corr/abs-1901-04098},
\begin{align*}
    \begin{bmatrix}
      u' \\
      v'
    \end{bmatrix}
    &=
    \begin{bmatrix}
      v \\
      -u - 8.53(u^2-1)v+1.2\sin\pl\dfrac{\pi}{5}t\pr
    \end{bmatrix},\qquad
    \begin{bmatrix}
      u(0)\\
      v(0)
    \end{bmatrix}
    =
    \begin{bmatrix}
      1.45 \\
      0
    \end{bmatrix}
\end{align*}
for $t\in[0,25]$.  As this IVP does not have an analytical solution, we compute reference solutions as described above.

We split the right-hand side function into linear and nonlinear slow and fast components, respectively,
\begin{align*}
    f^{\editmade{\mfs}} =
    \begin{bmatrix}
      v \\
      -u
    \end{bmatrix},\qquad
    f^{\editmade{\mff}} =
    \begin{bmatrix}
      0 \\
      - 8.53(u^2-1)v+1.2\sin(\frac{\pi}{5}t)
    \end{bmatrix}.
\end{align*}

\subsubsection{Pleiades}
\label{problem:pleiades}

The Pleiades problem is a special case of the general N-Body problem from \cite[Chapter~II.10]{hairer_solving_2009}, here comprised of seven bodies in two physical dimensions, resulting in fourteen position and fourteen velocity components ($p$ and $v$).  This problem has initial conditions
\begin{align*}
    p_1(0) &= \begin{bmatrix}
    3 & 3
    \end{bmatrix}, &
    p_2(0) &= \begin{bmatrix}
    3 & -3
    \end{bmatrix}, &
    p_3(0) &= \begin{bmatrix}
    -1 & 2
    \end{bmatrix}, &
    p_4(0) &= \begin{bmatrix}
    -3 & 0
    \end{bmatrix},\\
    p_5(0) &= \begin{bmatrix}
    2 & 0
    \end{bmatrix}, &
    p_6(0) &= \begin{bmatrix}
    -2 & 4
    \end{bmatrix}, &
    p_7(0) &= \begin{bmatrix}
    2 & 4
    \end{bmatrix} \\
    v_1(0) &= \begin{bmatrix}
    0 & 0
    \end{bmatrix}, &
    v_2(0) &= \begin{bmatrix}
    0 & 0
    \end{bmatrix}, &
    v_3(0) &= \begin{bmatrix}
    0 & 0
    \end{bmatrix}, &
    v_4(0) &= \begin{bmatrix}
    0 & -1.25
    \end{bmatrix},\\
    v_5(0) &= \begin{bmatrix}
    0 & 1
    \end{bmatrix}, &
    v_6(0) &= \begin{bmatrix}
    1.75 & 0
    \end{bmatrix}, &
    v_7(0) &= \begin{bmatrix}
    -1.5 & 0
    \end{bmatrix},
\end{align*}
and solutions were considered on the interval $t\in[0,3]$.

This IVP has no analytical solution, so we approximate reference solutions as described previously.  We split the right-hand side function component-wise, such that $f^{\editmade{\mfs}}$ contains the time derivatives of the positions, and $f^{\editmade{\mff}}$ contains the time derivatives of the velocities.

\subsubsection{FourBody3D}
\label{problem:fourbody3d}

The FourBody3D problem is another special case of the general N-Body problem, with four bodies in three spatial dimensions, defined in \cite{DBLP:journals/corr/abs-1901-04098}.  This problem has initial conditions
\begin{align*}
    p_1(0) &= \begin{bmatrix}
    0 & 0 & 0
    \end{bmatrix},\
    p_2(0) = \begin{bmatrix}
    4 & 3 & 1
    \end{bmatrix},\
    p_3(0) = \begin{bmatrix}
    3 & -4 & -2
    \end{bmatrix},\\
    p_4(0) &= \begin{bmatrix}
    3 & 4 & 5
    \end{bmatrix},\
    v_1(0) = v_2(0) = v_3(0) = v_4(0) = \begin{bmatrix}
        0 & 0 & 0
    \end{bmatrix},
\end{align*}
and solutions were considered on the interval $t\in[0,15]$.

This IVP has no analytical solution, and so reference solutions are approximated appropriately.  As with the Pleiades problem, we split the right-hand side function so that $f^{\editmade{\mfs}}$ contains the time derivatives of the positions, while $f^{\editmade{\mff}}$ contains the time derivatives of the velocities.

\subsection{Testing Suite}
\label{sec:testing-suite}

We evaluated the performance of each of the five proposed fast error estimation strategies, \emph{FS}, \emph{SA-mean}, \emph{SA-max}, \emph{LASA-mean}, and \emph{LASA-max}, over the above test problems, over the tolerances $\{10^{-3},10^{-5},10^{-7}\}$, and using each of the Constant-Constant, Linear-Linear, PIMR, and PIDMR controllers. \editmade{For simplicity of presentation, we always set $\tol^\mfs=\tol^\mff=\tfrac12\tol$}. When using controllers with extended histories, we used the Constant-Constant controller until a sufficient history had built up. Because we consider integer values of $M$, we take the \emph{ceil} of the value from the $M$ update functions \editmade{resulting from each controller}.

MRI-GARK \cite{sandu_class_2019} is the only family of infinitesimal methods we are aware of that has available embeddings for temporal error estimation. For our testing set we used:
\begin{itemize}
    \item MRI-GARK-ERK33, a four-stage third-order MRI-GARK method with a second-order embedding, which is explicit at each slow stage.
    \item MRI-GARK-ERK45a, a six-stage fourth-order MRI-GARK method with a third-order embedding, which is explicit at each slow stage.
    \item MRI-GARK-IRK21a, a four-stage second-order MRI-GARK method with a first-order embedding, which is explicit in three slow stages, and implicit in one.
    \item MRI-GARK-ESDIRK34a, a seven-stage third-order MRI-GARK method with a second-order embedding, which is explicit in four slow stages, and implicit in three.
\end{itemize}
We chose these methods because they cover a range of orders of accuracy, and include an equal number of explicit and implicit methods.
A post-publication correction was made to the method MRI-GARK-ERK45a \cite{sandu2018erk45afix}, where the embedding coefficient row of the $\Gamma^0$ matrix is replaced with
\begin{align*}
    \widehat{\Gamma}^0_6 = \begin{bmatrix}
    -\frac{1482837}{759520} & \frac{175781}{71205} & -\frac{790577}{1139280} & -\frac{6379}{56964} & \frac{47}{96} & 0
    \end{bmatrix}.
\end{align*}
We used the corrected coefficients in our tests.

In our numerical tests, we pair each MRI-GARK method with a fast explicit Runge--Kutta method of the same order.  Specifically, for MRI-GARK-IRK21a we use the second-order Heun--Euler method, for MRI-GARK-ERK33 and MRI-GARK-ESDIR34a we use the third-order Bogacki--Shampine method \cite{bogacki_32_1989}, and for MRI-GARK-ERK45a we use the \editmade{fourth-order Zonneveld method \cite{zonneveld1963automatic}}.

\editblock{
\subsection{Assessing Adaptive Performance}
\label{sec:assessing-adaptive-performance}

To assess performance of adaptivity controllers, we defined a brute force algorithm which, for a given IVP problem, numerical method, and $\tol$, finds the largest value of $H_n$ and smallest value of $M_n$ at each step that will result in an approximate solution with error at or very near $\tol$. We consider the resulting total numbers of $f^{\mfs}$ and $f^{\mff}$ evaluations, denoted $f^{\mfs}_{\text{opt}}$ and $f^{\mff}_{\text{opt}}$, to be the estimated optimal costs of an adaptive solve.  We provide a detailed description if this brute force algorithm (including pseudocode) in Appendix \ref{appendix:optimal-performance}.
}

\subsection{Controller Parameter Optimization}
\label{sec:controller-parameters-optimization}

Each of the controllers derived in Section \ref{sec:controllers} depend on a set of two to six free parameters.  To compare the ``best case'' for each controller and error measurement strategy, we first numerically optimized the controller parameters across our testing suite of problems, methods, tolerances, and fast error measurement strategies.

To measure the quality of a given set of parameters, we computed \editmade{three} metrics.  We first define the ``Error Deviation'' arising from a given adaptive controller on a given test $\tau$ as
\begin{equation}
\label{eq:relativeerrdeviation}
\textnormal{(Error Deviation)}_{\tau} =\log_{10}\left(\frac{\varepsilon}{\tol}\right),
\end{equation}
where $\varepsilon$ is defined as the maximum \editmade{relative} error over ten equally spaced points in the test problem's time interval, \editmade{which provides a robust error performance measurement across the entire interval.}
Here, a method that achieves solution accuracy close to the target $\tol$ will have (Error Deviation)$_{\tau}$ close to \editmade{zero. A value of (Error Deviation)$_{\tau}>0$ indicates $\varepsilon>\tol$ and (Error Deviation)$_{\tau}<0$ indicates $\varepsilon<\tol$.}

\editmade{Next, we use the total number of $f^{\mfs}$ and $f^{\mff}$ evaluations, $f^{\mfs}_{\text{evals}}$ and $f^{\mff}_{\text{evals}}$ (referred to as the ``slow cost'' and ``fast cost'', respectively), in relative factors to measure performance on a given test $\tau$.}
\begin{align}
  \textnormal{(\editmade{Slow} Cost Deviation)}_{\tau}=\frac{f^{\editmade{\mfs}}_{\textnormal{evals}}}{f^{\editmade{\mfs}}_{\textnormal{opt}}}, \label{eq:slow-cost} \\
  \textnormal{(\editmade{Fast} Cost Deviation)}_{\tau}=\frac{f^{\editmade{\mff}}_{\textnormal{evals}}}
  {f^{\editmade{\mff}}_{\textnormal{opt}}}, \label{eq:fast-cost}
\end{align}
A method that achieves \editmade{costs} close to optimal  will have \editmade{(Slow Cost Deviation)$_{\tau}$ and (Fast Cost Deviation)$_{\tau}$} close to one, although typically these values are significantly larger. \editmade{The primary purpose of this relative-to-optimal definition of cost is to accommodate the fact that different problems may require vastly different amounts of function evaluations, and we would like to consider the controllers' average performance across the testing suite.}

\editmade{To severely penalize parameter values that led to a lack of controller robustness}, we defined an optimization objective function to be
\begin{equation}
\label{eq:objective_function}
\begin{aligned}
    E(k) &= \sum_{\tau\in \text{test set}} E_{\tau}(k),\\
    E_{\tau}(k) &= \begin{cases}
    \editmade{10(\text{Slow Cost Deviation})_{\tau} + (\text{Fast Cost Deviation})_{\tau}} \\
    \ \ \ \ \editmade{+10(\text{Error Deviation})^2_{\tau}}, & \text{if $\tau$ finished}\\
    10^{10}, & \text{if $\tau$ failed}.
    \end{cases}
\end{aligned}
\end{equation}
\editmade{Here,} the contribution to the objective function for a particular test $\tau$ is small if the computed error is close to the target tolerance, and if the computational \editmade{costs are} not much larger than the ``optimal'' \editmade{values}. \editmade{We chose a factor of 10 for the Slow Cost Deviation to ensure that $f^{\mfs}_{\text{evals}}$ has a greater weight than $f^{\mff}_{\text{evals}}$ in the optimizer, and for the squared Error Deviation to prioritize computations that achieve the target accuracy.}

Due to the lack of differentiability in $E(k)$ due to failed solves and integer-valued $M$, we performed a simple optimization strategy consisting of an iterative search over the two- to six-dimensional parameter space.
We performed successive mesh refinement over an $n$-dimensional mesh $[0,1]^n$, with an initial mesh using a spacing of 0.2.
After evaluating the controller's performance on all sets of the parameters in the initial mesh, we refined the mesh around the parameter point having smallest objective function value with a mesh width of 0.4 in all $n$ directions, and a spacing of 0.04.
After evaluating the controller's performance on all of the parameters in this refined mesh, we refined the mesh a final time around the parameter point having smallest objective function value, with a mesh width of 0.08 in all $n$ directions and a spacing of 0.02.

\subsection{Fast error estimation strategy performance}
\label{sec:esf-strategies}

\begin{figure}[htb]
\centering
\includegraphics[width=0.75\textwidth]{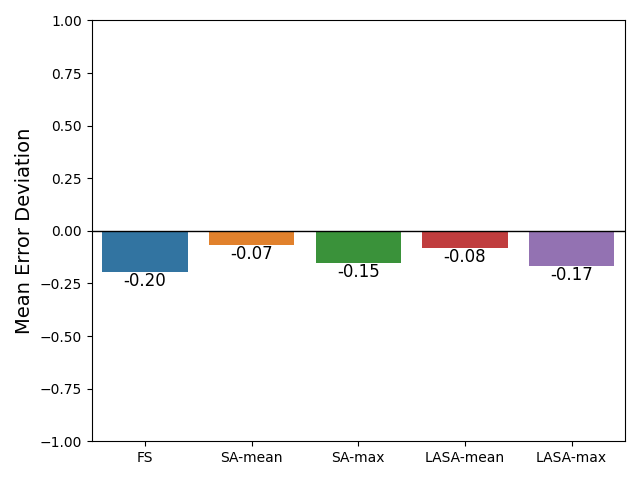}
\caption{\editmade{Mean Error Deviation arising from each fast error measurement strategy.  All proposed methods obtained results that achieved the desired tolerance.}}
\label{fig:errdeviation-by-esf}
\end{figure}

Our primary question for the quality of each of our fast error estimation strategies is how well it can estimate the solution error arising from approximation of each fast IVP.  Thus for each fast error measurement strategy, we define the average Error Deviation from the target tolerance as the average value of (Error Deviation)$_{\tau}$ over $\tau$ in a test set comprised of all combinations of our seven test problems, four IVP methods, three tolerances, and four controllers.  We plot these results in Figure \ref{fig:errdeviation-by-esf}, where we see that although each strategy followed a drastically different approach for error estimation, all \editmade{were able to achieve approximations that achieved the target solution accuracy.  However, we note that the FS, SA-max and LASA-max appear to have over estimated the fast error, leading to overly accurate results.}

\begin{figure}[htb]
\centering
\begin{subfigure}[b]{0.45\textwidth}
    \centering
    \includegraphics[width=\textwidth]{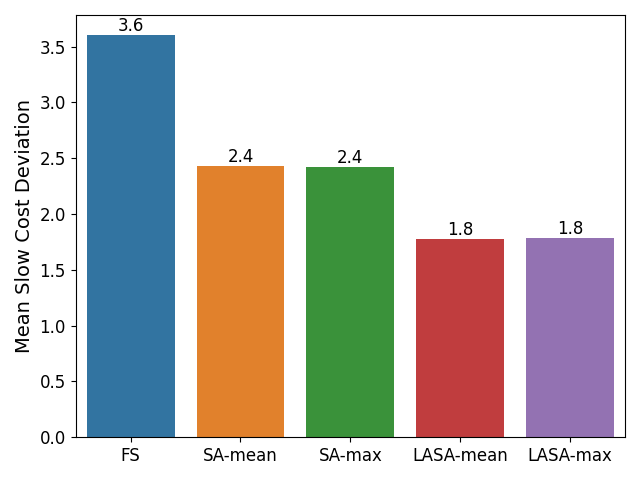}
    \caption{Slow Cost Deviation}
    \label{fig:esf_slow_cost}
\end{subfigure}
\hfill
\begin{subfigure}[b]{0.45\textwidth}
    \centering
    \includegraphics[width=\textwidth]{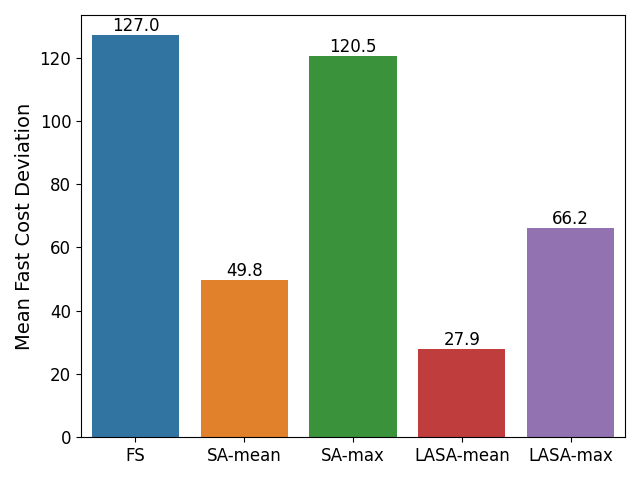}
    \caption{Fast Cost Deviation}
    \label{fig:esf_fast_cost}
\end{subfigure}
\caption{\editmade{Mean Slow and Fast Cost Deviation over test suite by each fast error measurement strategy.}}
\label{fig:cost-by-esf}
\end{figure}

Given that each fast error estimation strategy is able to obtain results of desirable accuracy, our second question focuses on the efficiency of using each approach in practice.  We thus define the relative cost of each fast error estimation strategy as the average value of both \editmade{(Slow Cost Deviation)$_{\tau}$ and (Fast Cost Deviation)$_{\tau}$} over $\tau$ in the same test set comprising all combinations of our seven test problems, four IVP methods, three tolerances, and four controllers.  We provide these plots in Figure \ref{fig:cost-by-esf}, where we see that all of the fast error estimation strategies had average \editmade{Slow Cost Deviation} within a factor of two from one another, with the LASA strategies providing the closest-to-optimal slow cost. Additionally, the ``LASA-mean'' strategy provided by far the closest-to-optimal fast cost by a significant margin.  This was expected, since the two LASA strategies were designed to minimize computational cost, yet their error estimates were sufficiently accurate.  \editmade{We believe that LASA-mean outperformed LASA-max because it provided a sharper estimate of fast solution error, as seen in Figure \ref{fig:errdeviation-by-esf}.  Based on these results}, in all subsequent numerical results we restrict our attention to the LASA-mean strategy alone.

\subsection{Optimized Controller Parameters}
\label{sec:controller-parameters}

After refining our focus to only the LASA-mean fast error estimation strategy, we re-ran the optimization approach described in Section \ref{sec:controller-parameters-optimization} to determine an ``optimal'' set of parameters for each $H$-$M$ controller.  These final parameters were:
\begin{itemize}
    \item Constant-Constant controller \eqref{eq:constant_constant_controller}:
    \begin{equation}
    \label{CCparams}
        k_1 = \editmade{0.42},\qquad k_2 = \editmade{0.44}.
    \end{equation}
    \item Linear-Linear controller \eqref{eq:linear_linear_controller}:
    \begin{equation}
    \label{LLparams}
        K_1 =
        \begin{bmatrix}
        \editmade{0.82} & \editmade{0.54}
        \end{bmatrix}^T,\qquad K_2 =
        \begin{bmatrix}
        \editmade{0.94} & \editmade{0.9}
        \end{bmatrix}^T.
    \end{equation}
    \item PIMR controller \eqref{eq:pimr_controller}:
    \begin{equation}
    \label{PIMRparams}
        K_1 =
        \begin{bmatrix}
        \editmade{0.18} & \editmade{0.86}
        \end{bmatrix}^T,\qquad K_2 =
        \begin{bmatrix}
        \editmade{0.34} & \editmade{0.80}
        \end{bmatrix}^T.
    \end{equation}
    \item PIDMR controller \eqref{eq:pidmr_controller}:
    \begin{equation}
    \label{PIDMRparams}
        K_1 =
        \begin{bmatrix}
        \editmade{0.34} & \editmade{0.10} & \editmade{0.78}
        \end{bmatrix}^T,\qquad K_2 =
        \begin{bmatrix}
        \editmade{0.46} & \editmade{0.42} & \editmade{0.74}
        \end{bmatrix}^T.
    \end{equation}
\end{itemize}

\subsection{Controller performance}
\label{sec:controller-performance}

With our chosen fast error estimation strategy and re-optimized parameters in place, we now \editmade{compare the performance of our newly-proposed $H$-$M$ controllers against the standard single-rate I, PI, and PID controllers, as well as Gustafsson's controller}.  For these tests, we utilized a subset of the testing suite above -- namely, for each controller we considered all seven test problems, all four IVP methods, and all three accuracy tolerances.  \editmade{For the single-rate controllers, we held $M=10$ constant for each test and used $\varepsilon^s$ as the temporal error estimate.}

\begin{figure}[htb]
\centering
\includegraphics[width=0.75\textwidth]{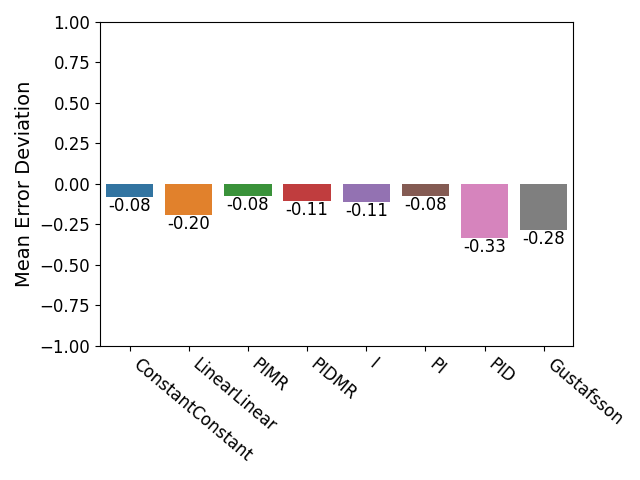}
\caption{Mean Error Deviation over test suite for each controller.}
\label{fig:errdeviation-by-controller}
\end{figure}

In Figure \ref{fig:errdeviation-by-controller}, we plot the average Error Deviation \eqref{eq:relativeerrdeviation} for each controller.  We see that all controllers again lead to solutions with Error Deviation \editmade{below zero}, implying the achieved errors for all approaches \editmade{achieve} their target tolerances \editmade{with a slight bias toward over-solving the problem.  However, we note that although the differences are small, the Linear-Linear controller provides solutions with errors furthest from $\tol$ \emph{of all of our proposed controllers}, with a mean Error Deviation of -0.20, indicating that the solution had error approximately $6\times10^{-4}$ when $\tol=10^{-3}$, or $6\times10^{-8}$ when $\tol=10^{-7}$, which are still well within range of the target tolerance. We note that both the PID and Gustafsson single-rate controllers achieved solutions that were considerably more accurate than requested, although even those were within a reasonable range of the tolerance.}

\begin{figure}[htb]
\centering
\begin{subfigure}[b]{0.45\textwidth}
    \centering
    \includegraphics[width=\textwidth]{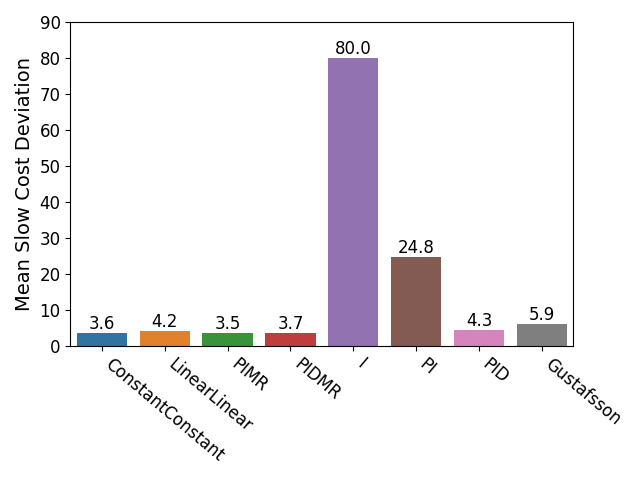}
    \caption{Slow Cost Deviation}
    \label{fig:controller_slow_cost}
\end{subfigure}
\hfill
\begin{subfigure}[b]{0.45\textwidth}
    \centering
    \includegraphics[width=\textwidth]{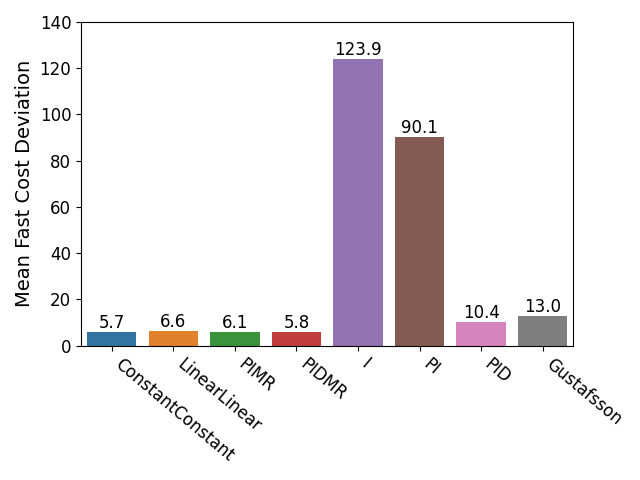}
    \caption{Fast Cost Deviation}
    \label{fig:controller_fast_cost}
\end{subfigure}
\caption{Mean Slow and Fast Cost Deviation over test suite by each controller.}
\label{fig:cost-by-controller}
\end{figure}

\editmade{In Figure \ref{fig:cost-by-controller}, we plot the average cost deviation for each controller, over all of the combinations of test problems, methods, and tolerances. The proposed multirate controllers demonstrated comparable computational cost across our ODE test suite, with differences between the best and worst controllers of only 20\% in terms of the average Slow Cost Deviation and only 16\% in terms of the average Fast Cost Deviation.
Meanwhile, the proposed multirate controllers uniformly outperformed their single-rate counterparts. Of those, the PID controller had the best performance which was on-par in average Slow Cost Deviation  but 58\% worse in average Fast Cost Deviation than the worst-performing multirate controller (Linear-Linear). Gustafsson's controller had a slightly worse performance, and the I and PI controller's were dramatically worse. Additionally, we note that the I controller failed all runs on the Forced Van der Pol problem when using MRIGARKIRK21a and MRIGARKESDIRK34a, and so we did not include those runs in these I controller averages.}

\editmade{Based on these results, it is clear that the proposed multirate controllers all show excellent performance, with no single method outperforming another.  Thus, for relatively simple IVPs we recommend the Constant-Constant controller due to its simplicity, whereas for more complex IVPs we recommend testing with each multirate controller.}

\subsection{Multirate controller performance deep dive}
\label{sec:controller-comparisons}

\editblock{
  The previous results focused on averaged controller performance over a wide range of problems on which the controller parameters had already been optimized.  In this section we instead compare the performance of our proposed controllers on a new test problem for which our adaptivity controllers have not been optimized, and that should thoroughly exercise their ability to adapt step sizes at both the fast and slow time scales.  Thus this should provide an unbiased challenge problem on which we may compare controller performance, while also allowing a deeper dive into controller behavior.

  We adapt the stiff Brusselator example from Section \ref{problem:brusselator} to a 1D reaction-diffusion setting with time-varying coefficients,
  \begin{align*}
    \partial_t u &= d(t)\,\partial_{xx} u + r(t)\left(a - (w+1)u + u^2v\right),\\
    \partial_t v &= d(t)\,\partial_{xx} v + r(t)\left(uw-u^2v\right),\\
    \partial_t w &= d(t)\,\partial_{xx} w + r(t)\left(\frac{b-w}{\epsilon} - uw\right),
  \end{align*}
  for $(t,x) \in (0,2) \times (0,1)$, with initial conditions
  \begin{equation*}
    u(0) = 1.2 + 0.1\sin(\pi x), \; v(0) = 3.1 + 0.1\sin(\pi x), \; w(0) = 3 + 0.1\sin(\pi x),
  \end{equation*}
  stationary boundary conditions,
  \begin{equation*}
    \partial_t u(t,0) = \partial_t u(t,1) = \partial_t v(t,0) = \partial_t v(t,1) = \partial_t w(t,0) = \partial_t w(t,1) = 0,
  \end{equation*}
  time-varying coefficient functions
  \begin{equation*}
    d(t) = 0.006 + 0.005 \cos(\pi t), \qquad r(t) = 0.6 + 0.5\cos(4\pi t),
  \end{equation*}
   and parameters $a=1$, $b=3.5$.  Here, we increase the stiffness of the problem by setting $\varepsilon=0.001$ (previously this was 0.01).  We partition this problem such that $f^{\mfs}$ corresponds to the diffusion terms, while $f^{\mff}$ corresponds to the reaction terms.  We note that for the above values, the diffusion coefficient $d(t)$ lies within $(0.001, 0.011)$ and the reaction coefficient $r(t)$ lies within $(0.1, 1.1)$, but that the frequencies of these oscillations differ, leading to coefficient ratios $d(t)/r(t)$ that range from approximately $10^{-3}$ to $10^{-1}$.  We thus expect that each of our adaptivity controllers will need to vary both $H$ and $M$ to accurately track the multirate solutions.

   In lieu of averaging performance values across a multitude of methods, we focus on only MRIGARKERK45a here, although we note that the results are similar when using other multirate methods.
}

\begin{figure}[htb]
\centering
\begin{subfigure}[b]{0.45\textwidth}
    \centering
    \includegraphics[width=\textwidth]{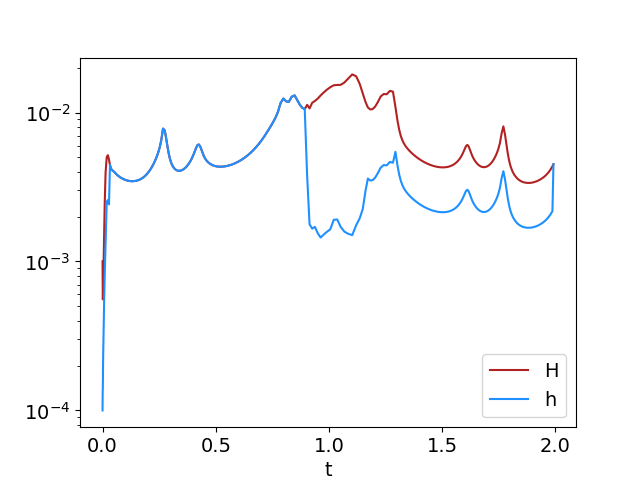}
    \caption{Constant-Constant}
    \label{fig:constantconstant-H-h}
\end{subfigure}
\hfill
\begin{subfigure}[b]{0.45\textwidth}
    \centering
    \includegraphics[width=\textwidth]{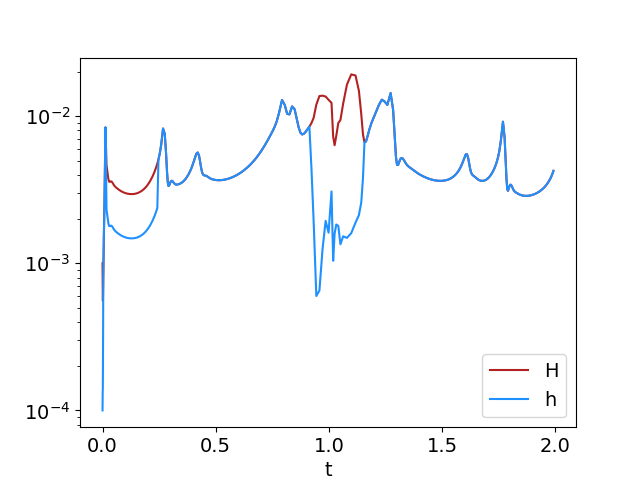}
    \caption{Linear-Linear}
    \label{fig:linearlinear-H-h}
\end{subfigure}
\newline
\begin{subfigure}[b]{0.45\textwidth}
    \centering
    \includegraphics[width=\textwidth]{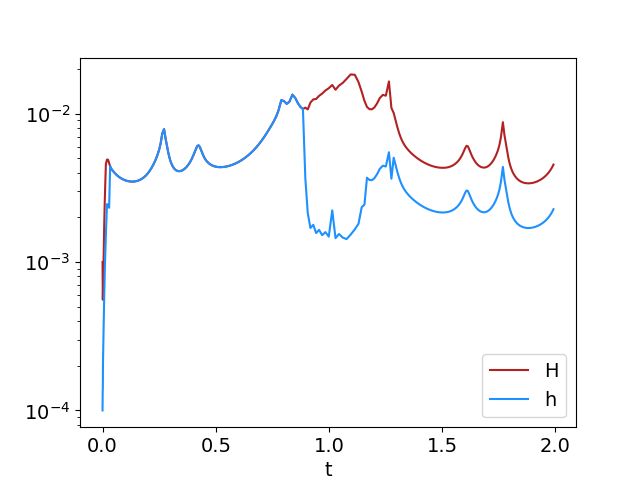}
    \caption{PIMR}
    \label{fig:pimr-H-h}
\end{subfigure}
\hfill
\begin{subfigure}[b]{0.45\textwidth}
    \centering
    \includegraphics[width=\textwidth]{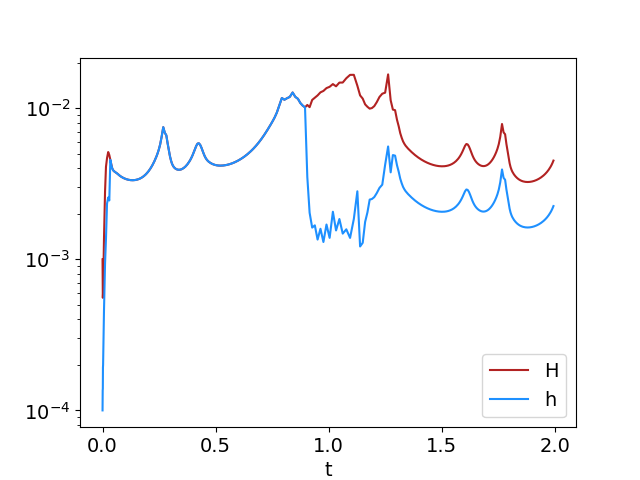}
    \caption{PIDMR}
    \label{fig:pidmr-H-h}
\end{subfigure}
\caption{\editmade{$H_n$ and $h_n$ over time for each multirate controller with a tolerance of $10^{-4}$ on the 1D stiff Brusselator problem.}}
\label{fig:brusselator1d-H-h}
\end{figure}

\editmade{In Figure \ref{fig:brusselator1d-H-h} we plot the time step sizes $H_n$ and $h_n$ over time for each of our controllers with a tolerance of $10^{-4}$.  We can see that for the first half of each simulation the problem did not exhibit multirate behavior, so all controllers varied $H$ similarly and set $M=1$ (except for Linear-Linear, that showed a brief initial period with $M>1$). At approximately $t=0.9$, some stiffness arises in the reaction network and the controllers all respond by increasing $M$ and thus decreasing $h$. We can see the effect of some failed steps at $t=1.0$ in the Linear-Linear controller, where it decreases $H$ and rapidly adjusts the value of $M$, while the other controllers more smoothly adjust $H$ with some higher-frequency changes to $M$. Once the period of stiffness ends (around $t=1.25$), the Linear-Linear controller resets $M$ to 1, while the other controllers maintain a small value of $M=2$.}
\begin{figure}[htb]
\centering
\begin{subfigure}[b]{0.45\textwidth}
    \centering
    \includegraphics[width=\textwidth]{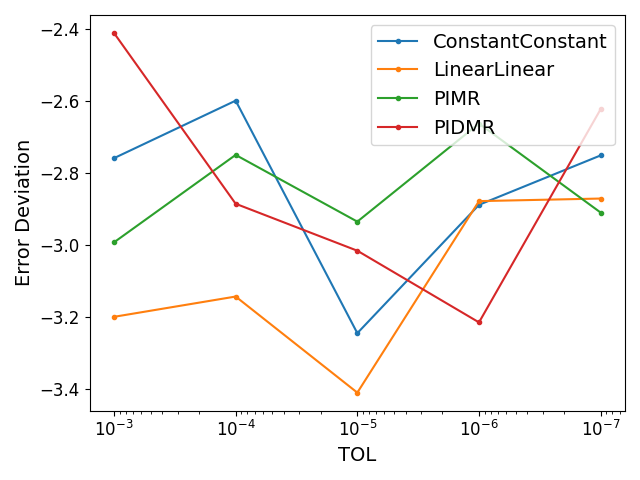}
    \caption{}
    \label{fig:brusselator1d-errdev}
\end{subfigure}
\newline
\begin{subfigure}[b]{0.45\textwidth}
    \centering
    \includegraphics[width=\textwidth]{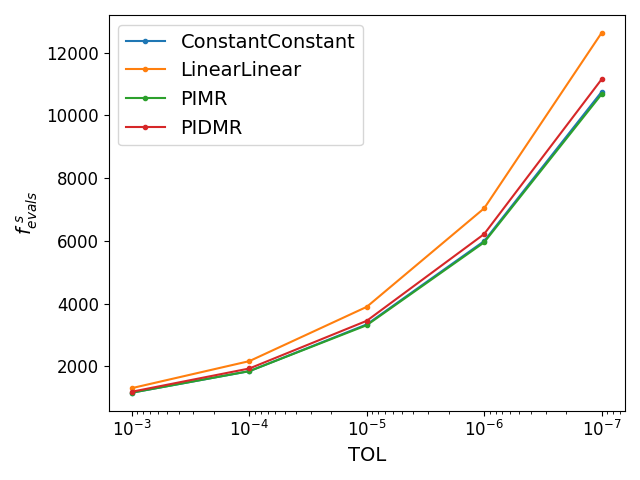}
    \caption{}
    \label{fig:brusselator1d-slow-cost}
\end{subfigure}
\hfill
\begin{subfigure}[b]{0.45\textwidth}
    \centering
    \includegraphics[width=\textwidth]{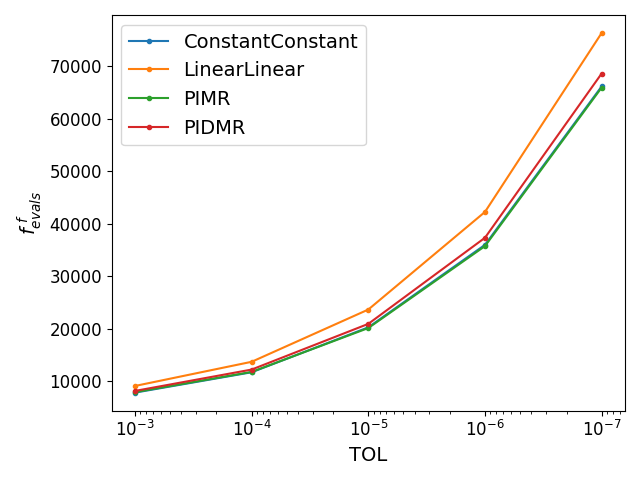}
    \caption{}
    \label{fig:brusselator1d-fast-cost}
\end{subfigure}
\caption{\editmade{(a) Error Deviation, (b) Slow function evaluations, and (c) Fast function evaluations vs. $\tol$ for each multirate controller on the 1D stiff Brusselator problem.}}
\label{fig:brusselator1d-stats}
\end{figure}

\editblock{
In Figure \ref{fig:brusselator1d-stats} we examine the performance of each controller as the tolerance is varied, with plots of the Error Deviation, the total slow function evaluations, and the total fast function evaluations for each controller.  All multirate controllers over-solved the problem, providing solutions two to four orders of magnitude more accurate than the chosen $\tol$, however there was no consistent pattern as to which controller performed the best for loose/tight $\tol$ values.  Similarly, the controllers show comparable performance in terms of computational cost. Only the Linear-Linear controller had noticeably more function evaluations than the others, with approximately 13\% more slow and 11\% more fast function evaluations than the other controllers.
}


\section{Conclusions}
\label{sec:conclusions}

We followed the technique of Gustafsson \cite{gustafsson_control-theoretic_1994} to develop controllers that approximate the fast and slow principal error functions for \editmade{multirate infinitesimal} methods.  To this end, we developed piecewise constant and linear approximations for each principal error function.  We then combined these approximations using pairs of piecewise polynomial approximations to the principal error functions with like degree to construct Constant-Constant and Linear-Linear controllers for both the slow time step size $H$ and the multirate ratio $M$, within \editmade{multirate infinitesimal} methods.

To assess the reliability and to measure the efficiency of these proposed controllers, we devised a large testing suite encompassing seven multirate test problems, four MRI methods and three accuracy tolerances.  In order to measure method efficiency, we developed an algorithm to determine the best-case pair of $H_n$ and $M_n$ values for each testing combination.

In our initial tests, however, we found that controllers with polynomial approximations of degree two or larger to the principal error function tended to constrain step size changes too tightly, leading to a large number of solver failures. To address this issue, we introduced the PIMR controller, formed by taking the Linear-Linear controller and removing dependence on $H$ and $M$ terms \editmade{prior to $H_n$ and $M_n$}, and the PIDMR controller, an extension to the PIMR controller with an increased error history. These controllers were developed to have similar structures to the existing single-rate PI and PID controllers. Through these modifications, the PIMR and PIDMR controllers can react more swiftly to a problem's influences on the multirate step size(s).

While estimation of the slow error $\varepsilon^{\editmade{\mfs}}$ is straightforward for multirate methods with embeddings, we developed multiple strategies to estimate the fast error value $\varepsilon^{\editmade{\mff}}$, including the Full-Step, Stage-Aggregate, and Local-Accumulation-Stage-Aggregate strategies.
These strategies trade off differing levels of computational effort with the expected accuracy in their estimation.  However, when examining the performance of these approaches over our testing suite, we found that the Local-Accumulation-Stage-Aggregate strategy with mean aggregation offered an ideal combination of low cost and reliable accuracy to the target tolerance, and we therefore recommend it for practitioners.

We then evaluated the performance of each controller over our testing set, finding that the Constant-Constant, Linear-Linear, PIMR, and PIDMR controllers each tended to achieve solutions close to the target tolerance.  \editmade{While our proposed controllers perform similarly, they greatly outperform existing single-rate I, PI, Gustafsson's controllers, and slightly outperform the PID controller, in terms of both slow and fast function evaluations on average.}

\editmade{Finally, we evaluated our controllers on a stiffer, PDE version of the Brusselator problem. We saw that the controllers adjust both $H$ and $M$ in response to a period of increased stiffness, and readjusted after the period had ended. The controllers had a roughly equivalent computational cost in solving this problem, with the Linear-Linear controller consistently experiencing a slightly higher cost. Each controller tended to over-solve the problem, giving solutions with errors two to four orders of magnitude lower than the chosen value of $\tol$.}

Significant work remains in the area of temporal adaptivity for multirate methods.  MRI-GARK is the only infinitesimal method family we have found that includes embeddings. Thus, embedded methods from other \editmade{multirate infinitesimal} families need to be derived, along with a greater ecosystem of embedded MRI-GARK methods that focus more specifically on performance in a temporally adaptive context.

We additionally note that controllers which update $H$ and $M$ for each slow multirate step may not be the most efficient choice.  We focused on $H$-$M$ controllers, as those give $H$ values that are an integer multiple of $h$ values and result in a simpler method implementation; however, controllers may be created instead for $H$ and $h$ by following the steps outlined in this paper, replacing $\frac{H}{M}$ with $h$ in the early steps.  Perhaps the increased flexibility arising from real-valued $h$ could lead to efficiency improvements over the integer-valued $M$ approaches here.


\section{Acknowledgments}
\label{sec:acknowledgments}

We would like to thank the Virginia Tech Computational Science Laboratory for the test problem repository \cite{DBLP:journals/corr/abs-1901-04098}, which proved to be a huge help in finding test problems to evaluate our controllers.
We would like to thank Arash Sarshar specifically for his help in understanding and debugging our implementations of some MrGARK and MRI-GARK methods.
We would also like to thank Steven Roberts for his help in reviewing this manuscript and his help in understanding MRI-GARK methods.


\bibliographystyle{siamplain}
\bibliography{references}
\newpage
\appendix
\section{Optimal Performance Estimation Algorithms}
\label{appendix:optimal-performance}

In order to compare the performance of our proposed adaptive controllers and error estimation algorithms, we create a baseline set of ``optimal'' cost values. \editmade{We note that in an practice the most} computationally efficient values of $H_n$ and $M_n$ will depend on a number of factors, including: the IVP itself, the multirate method under consideration, \editmade{the cost of any implicit solvers at either time scale,} the desired solution accuracy, and even the relative computational cost of evaluating the slow and fast right-hand side functions, $f^{\editmade{\mfs}}$ and $f^{\editmade{\mff}}$.  Furthermore, even these optimal values of $H_n$ and $M_n$ will vary as functions of time throughout the simulation, particularly for nontrivial multirate problems.

\editmade{In this work, we} define the optimal \editmade{cost as the minimal number of $f^{\mfs}$ and $f^{\mff}$ evaluations required to reach the end of the time interval, where} each step results in local error estimates that achieve the chosen tolerance, and with each step locally optimal with respect to a prescribed computational efficiency measurement.  For the sake of simplicity, we define this efficiency measurement as
\begin{align}
  \label{eq:efficiency}
  \textnormal{efficiency} &= \frac{H_n}{\textnormal{cost}}, \\
  \label{eq:cost}
  \textnormal{cost}  &=\textnormal{slowWeight}\cdot f^{\editmade{\mfs}}_{\textnormal{evals}} + f^{\editmade{\mff}}_{\textnormal{evals}},
\end{align}
where $f^{\editmade{\mfs}}_{\textnormal{evals}}$ and $f^{\editmade{\mff}}_{\textnormal{evals}}$  are the total number of $f^{\editmade{\mfs}}$ and $f^{\editmade{\mff}}$ evaluations for the multirate time step, respectively.  Here, ``slowWeight'' provides a problem-specific factor that encodes the relative costs of $f^{\editmade{\mfs}}$ and $f^{\editmade{\mff}}$.  We note that for any given simulation, this value could itself depend on numerous un-modeled factors, such as the IVP under consideration, its numerical implementation, and even the computational hardware.  However, irrespective of the ``slowWeight'' value used, for a given slow step size $H_n$, a method that results in a smaller overall ``cost'' corresponds with increased efficiency.

We chose the definitions \eqref{eq:efficiency}-\eqref{eq:cost} because, in the goal of achieving the cheapest possible solve of a given IVP to a given tolerance, we want as large of step sizes as possible, and as small of costs as possible.
If a step is rather expensive, e.g. if the value of $M_n$ is high, the step can still achieve a high efficiency if the step size was large.
Eventually, once the errors arising from the fast time scale are sufficiently small for a given method or problem, additional increases to $M_n$ will not improve the overall accuracy and will thereby lead to decreased efficiency.  Similarly, $H_n$ will be bounded from above due to accuracy considerations, and although decreasing $H_n$ below this bound may allow for a smaller $M_n$, the overall efficiency could decrease.


With these definitions in place, our approach to find the optimal set of $H$-$M$ pairs is shown \editmade{in pseudocode representation} in Algorithms \ref{alg:optimalHM} and \ref{alg:findH}.  This is essentially a brute-force mechanism to rigorously determine the best-case values for multirate adaptivity algorithms.
The function ``ComputeReferenceSolution'' is a black box that computes the reference solution at a desired time $t_i+H$ and is assumed to be more accurate than the ``ComputeStep'' function.
The function ``ComputeStep'' is a black box function that takes one step with the given method from $t_i$ to $t_i+H$ and returns the total slow and fast function calls, the error in the step's solution, and the solution itself.
For a given IVP, initial condition, and initial time, the algorithm iterates over increasing values of the integer multirate ratio $M_n$ and uses the given multirate method to find the maximal step size $H_n$ for each $M_n$ which gives an error close to the chosen tolerance via a binary search process, stopping when the interval width is smaller than a relative tolerance $H_{tol}$ of the midpoint of the interval.
Once the efficiency from increasing $M_n$ decreases below some relative tolerance $\textnormal{eff}_{rtol}$ of the maximum so far found, the solution is moved forward based on the most efficient $(H_n, M_n)$ pair and repeats, iterating until the algorithm reaches the end of the given time window $[t_0,t_f]$.

Algorithms \ref{alg:optimalHM} and \ref{alg:findH} are rather costly and the results from one run are specific to the IVP, method, and other parameters. For consistency, we always run the algorithm with the parameters $\textnormal{slowWeight}=10$, $H_{fine}=10^{-10}$, $H_{tol}=10^{-5}$, $H_{interval}=10^{-1}$, $M_{max\_iter}=400$, $M_{min\_iter}=10$, $\textnormal{eff}_{rtol}=10^{-1}$\emph{, and a sixth order explicit RK method \cite{verner1978} with small time steps for reference solutions.} We further note that the resulting ``optimal'' \editmade{total $f^{\mfs}$ and $f^{\mff}$ evaluations across each most efficient step} found by this algorithm achieve a cost that is nearly impossible for a time adaptivity controller to reach in practice, and should thus be considered a best possible scenario.

\begin{algorithm}[htbp]
\KwResult{Optimal $H$ array $H_{opt}$, Optimal $M$ array $M_{opt}$, Total $f^{\editmade{\mfs}}$ evaluations $f^{\editmade{\mfs}}_{opt}$, Total $f^{\editmade{\mff}}$ evaluations $f^{\editmade{\mff}}_{opt}$.}
Given an IVP, multirate method, error tolerance $tol$, weight factor slowWeight, initial condition $y_0$, time interval $\{t_0,t_f\}$, minimum step $H_{fine}$, binary search $H$ tolerance $H_{tol}$, binary search $H$ interval width $H_{interval}$, $M$ maximum $M_{max\_iter}$, $M$ minimum $M_{min\_iter}$, and relative efficiency tolerance $\eff_{rtol}$:\;
\vskip 0.25cm

$f^{\editmade{\mfs}}_{opt}\gets0$, $f^{\editmade{\mff}}_{opt}\gets0$, $i\gets 0$, $t\gets t_0$, $y_i\gets y_0$\;

\While{$t+H_{fine}<t_f$}{

	empty $H_{array}$, $M_{array}$, $\textnormal{eff}_{array}$, $f^{\editmade{\mfs}}_{evals,array}$, $f^{\editmade{\mff}}_{evals,array}$, $y_{array}$\;

	$M\gets1$\;

	\While{$M<M_{max\_iter}$}{
		$H$, $\eff$, $f^{\editmade{\mfs}}_{evals}$, $f^{\editmade{\mff}}_{evals}$, $y_{temp}\gets$ FindH(IVP, method, $tol$, slowWeight, $y$, $t$,\\
		\phantom {........................................................} $H_{fine}$, $M_{new}$, $H_{tol}$, $H_{interval}$)\;

		\eIf{$\frac{\textnormal{eff$_{array}$.max() - eff}}{\textnormal{eff$_{array}$.max()}} > \eff_{rtol}$ \textnormal{and} $M>M_{min\_iter}$}{
			break\;
		}{
			$M_{array}$.append($M$)\;

			$H_{array}$.append($H$)\;

			$\textnormal{eff}_{array}$.append(eff)\;

			$f^{\editmade{\mfs}}_{evals,array}$.append($f^{\editmade{\mfs}}_{evals}$)\;

			$f^{\editmade{\mff}}_{evals,array}$.append($f^{\editmade{\mff}}_{evals}$)\;

			$y_{array}$.append($y_{temp}$)
		}
	}
	$opt\_idx\gets$eff$_{array}$.indexOf(eff$_{array}$.max())\;

	$H_{opt}$.append($H_{array}[opt\_idx]$)\;

	$M_{opt}$.append($M_{array}[opt\_idx]$)\;

	$f^{\editmade{\mfs}}_{opt}\gets f^{\editmade{\mfs}}_{opt} + f^{\editmade{\mfs}}_{evals,array}[opt\_idx]$

	$f^{\editmade{\mff}}_{opt}\gets f^{\editmade{\mff}}_{opt} + f^{\editmade{\mff}}_{evals,array}[opt\_idx]$

	$t\gets t+H_{array}[opt\_idx]$\;

	$y\gets y_{array}[opt\_idx]$\;
}
\caption{Optimal H-M Search Algorithm}
\label{alg:optimalHM}
\end{algorithm}

\begin{algorithm}[htbp]
\KwResult{Maximal step size $H$ giving error close to $tol$, efficiency of computation $\eff$ using step size $H$, number of $f^{\editmade{\mfs}}$ evaluations $f^{\editmade{\mfs}}_{evals}$, number of $f^{\editmade{\mff}}$ evaluations $f^{\editmade{\mff}}_{evals}$, computed solution $y_{i+1}$ using step size $H$.}
Given an IVP, multirate method, error tolerance $tol$, weight factor slowWeight, initial condition $y_i$, initial time $t_i$, minimum step $H_{fine}$, multirate factor $M$, binary search $H$ tolerance $H_{tol}$, and binary search $H$ interval width $H_{interval}$:\;

$y_{ref}\gets$ ComputeReferenceSolution(IVP, $y_i$, $t_i$, $H_{fine}$)\;

$err$, $f^{\editmade{\mfs}}_{evals}$, $f^{\editmade{\mff}}_{evals}$, $y_{i+1}\gets$ ComputeStep(IVP, method, $y_i$, $t_i$, $H_{mid}$, $M$, $y_{ref}$)\;

$cost\gets\textnormal{slowWeight}\cdot f^{\editmade{\mfs}}_{evals}+f^{\editmade{\mff}}_{evals}$\;

eff $\gets H/\textnormal{cost}$\;

\eIf{$err < tol$}{

    $H_{right}\gets 0$\;

	\While{$err<tol$ \textnormal{and} $t_i+H_{right}<t_f$}{
		$H_{left} \gets H_{right}$\;

		$H_{right} \gets \min(H_{right}+H_{interval},\ t_f-t_i)$\;

		$H_{mid} \gets \frac{1}{2}(H_{left}+H_{right})$\;

		$y_{ref}\gets$ ComputeReferenceSolution(IVP, $y_i$, $t_i$, $H_{right}$)\;

		$err$, $f^{\editmade{\mfs}}_{evals}$, $f^{\editmade{\mff}}_{evals}$, $y_{i+1}\gets$ ComputeStep(IVP, method, $y_i$, $t_i$, $H_{right}$, $M$, $y_{ref}$)\;

		$cost\gets\textnormal{slowWeight}\cdot f^{\editmade{\mfs}}_{evals}+f^{\editmade{\mff}}_{evals}$\;

		eff $\gets H/\textnormal{cost}$\;

		$n\gets n+1$\;
	}
	\eIf{$err > tol$}{
    	\While{$(H_{right}-H_{left})/ H_{mid} > H_{tol}$}{
    		$H_{mid}\gets \frac{1}{2}(H_{left} + H_{right})$\;

    		$y_{ref}\gets $ComputeReferenceSolution(IVP, $y_i$, $t_i$, $H_{right}$)\;

    		$err,f^{\editmade{\mfs}}_{evals},f^{\editmade{\mff}}_{evals},y_{i+1}\gets$ ComputeStep(IVP,method,$y_i,t_i,H_{mid},M,y_{ref}$)\;

    		$cost\gets\textnormal{slowWeight}\cdot f^{\editmade{\mfs}}_{evals}+f^{\editmade{\mff}}_{evals}$\;

    		eff $\gets H/\textnormal{cost}$\;

    		\eIf{$err\le tol$}{
    			$H_{left}\gets H_{mid}$\;
    		}{
    			$H_{right}\gets H_{mid}$\;
    		}
    	}
    	$H\gets H_{left}$
	}{
	    $H\gets H_{right}$
	}
}{
	Failure ($H_{fine}$ was insufficiently small).\;
}
\caption{FindH Algorithm}
\label{alg:findH}
\end{algorithm}

\end{document}